\documentstyle[a4,leqno,amsfonts,12pt]{article}
\newcommand{\C}{{\bf C}}

\newcommand{\R}{{\bf R}}
\newcommand{\He}{{\bf H}}
\newcommand{\N}{{\bf N}}
\newcommand{\Z}{{\cal Z}}

\newcommand{\D}{{\bf D}}
\newcommand{\de}{\delta}
\newcommand{\De}{\Delta}
\newcommand{\mb}{\mbox}

\newcommand{\beq}{\begin{equation}}
\newcommand{\eeq}{\end{equation}}
\newcommand{\oge}{\succeq}
\newcommand{\ole}{\preceq}

\newcommand{\ov}{\overline}
\newcommand{\al}{\alpha}

\newcommand{\Om}{\Omega}
\newcommand{\om}{\omega}
\newcommand{\z}{\zeta}

\newcommand{\ga}{\gamma}
\newcommand{\la}{\lambda}
\newcommand{\Ga}{\Gamma}
\newcommand{\si}{\sigma}
\newcommand{\Si}{\Sigma}
\newtheorem{th}{Theorem}

\newtheorem{lem}{Lemma}

\newcommand{\ueberschrift}{\bigskip\goodbreak\noindent\bigskip}
\newcounter{theabsatz}
\newcommand{\absatz}[1]{\stepcounter{theabsatz} \ueberschrift
                           {\large \bf \arabic{theabsatz}. {#1}} \setcounter{equation}{0}}

\begin{document}
\mathsurround=2pt

\begin{center}
{\large \bf ON APPROXIMATION OF CONTINUOUS FUNCTIONS BY ENTIRE
FUNCTIONS OF EXPONENTIAL TYPE  ON SUBSETS OF
 THE REAL LINE}\\[3ex]

 Vladimir V.
Andrievskii$^*$
\end{center}

\begin{abstract}
We generalize the classical Bernstein theorem concerning the
 constructive description of classes
 of functions uniformly continuous on the real line.
 The approximation of continuous bounded functions  by entire functions
 of exponential type   on  unbounded closed proper subsets of
 the real line is studded.
\end{abstract}

\footnotetext{
 $^*$ The work    was supported in part by   NSF grant
DMS-0554344 and DFG grant BL 272/10-1.

Date received:$\quad$. Communicated by


 }

\absatz{Introduction and Main Results}

For a closed unbounded set $E\subset\R$, denote by $BC(E)$ the
class of (complex-valued) functions which are bounded and
continuous on $E$. Let $E_\si$ be the class of entire functions of
exponential type at most $\si>0$ and let
$$A_\si(f,E):=\inf_{g\in E_\si}||f-g||_{C(E)}\quad (f\in BC(E)),$$
 where $||\cdot||_{C(E)}$ means the uniform norm over $E$.

 The classical Bernstein direct and inverse theorems (see
 \cite[p.p. 257, 340]{tim}) describe the relations between the
 smoothness of $f\in BC(\R)$ and the rate of decrease of
 $A_\si(f,\R)$ as $\si\to\infty$. In particular,  from
  Bernstein's results it follows that for $f\in BC(\R)$ and $0<\al<1$,
 \beq\label{1.01}
 A_\si(f,\R)=O(\si^{-\al})\quad\mb{as }\si\to\infty
 \eeq
 if and only if
\beq\label{1.02} \om_{f,\R}(\de)=O(\de^\al)\quad \mb{as }\de\to
+0,
 \eeq
 where
 $$\om_{f,\R}(\de):=\sup_{x_1,x_2\in\R\atop|x_1-x_2|\le\de}|f(x_2
 )-f(x_1)|\quad (\de>0).
 $$
 The main objective of this paper is to extend  Bernstein's
 results to the case where the function is considered on  a proper subset
 of $\R$. We  mostly focus on the case where the number
 of components of $E$  is infinite.
 Some aspects of this problem are considered in the recent papers
 by Shirokov \cite{shi03}, \cite{shi04}.

 Let
 $$d(A,B):=\mb{dist}(A,B)=\inf_{z\in A,\z\in B}|z-\z|\quad
 (A,B\subset\C),
 $$
 and let $|B|$ denote the one-dimensional Lebesgue measure
 (length) of $B\subset\C$.

 Unless otherwise stated,  we denote by $
 C,C_1,C_2,\ldots$ positive constants
 that are either  absolute or depend  on $E$ only.

 The set $E^*:=\R\setminus E$ consists of a finite or infinite
 number  of disjoint open intervals $J_j=(a_j,b_j)$. In
 the reminder of this paper we assume that $E$ possesses the following
 two properties: for any $j$ under consideration,
 \beq\label{1.1}
 |J_j|\le C_1,
 \eeq
 \beq\label{1.2}
 \sum_{k\neq j}\left(\frac{|J_k|}{d(J_k,J_j)}\right)^2\le C_2.
 \eeq
 We use the following examples to illustrate the forthcoming
 results and constructions. The examples show that the number of
 ``holes" $J_j$ can be infinite.

\noindent  {\bf Example 1.} Let
  $$d_{l-1}<c_l<d_l<c_{l+1}\quad
  (l=0,\pm1,\pm2,\ldots)$$
  be such that
  $$d_l-c_l\ge C_3,\quad c_{l+1}-d_l\le C_4.
  $$
  Then, the set
  $$E_1=\cup_{l=-\infty}^\infty[c_l,d_l]$$
  satisfies (\ref{1.1}) and (\ref{1.2}).

\noindent {\bf Example 2.} A direct computation shows that
 the set $E_2=\R\setminus E_2^*$, where
 \begin{eqnarray*}
 E_2^*=\cup_{j=-\infty}^\infty\cup_{k=2}^\infty&&\{
 ( 2j+2^{-k}(1-k^{-1}),2j+2^{-k})\\
 &&\cup (2j-2^{-k},
 2j-2^{-k}(1-k^{-1}) )\}
 \end{eqnarray*}
 also satisfies (\ref{1.1}) and (\ref{1.2}).

 In the case of polynomial approximation of continuous functions
 on a finite interval $[a,b]\subset\R$, the special role of the
 endpoints $a$ and $b$ is well-known.
 An elegant idea,
  suggested in \cite{dittot}, is to
 introduce a new modulus of continuity by using the
 distance between the points on $[a,b]$ that is not Euclidean.
  In the case of
  entire function approximation on $E$ the endpoints of $J_j$
 also play a special role.
 We capture this effect by making use of a special distance
 between points of $E$
 in the definition of the modulus of continuity of a
 function $f\subset BC(E)$. This distance
  is defined as follows. Let $\He:=\{ z:\, \Im z >0\}$ be the upper half-plane.
 According to Levin \cite{lev89}
 there exist vertical intervals $J_j'=(u_j,u_j+iv_j],\,
 u_j\in\R,v_j>0$ and a conformal mapping
 $$\phi:\He\to\He_E:=\He\setminus(\cup_jJ_j')$$
 normalized by $\phi(\infty)=\infty,\phi(i)=i$ such
 that $\phi$ can be extended continuously to
 $\ov{\He}$ and it
 satisfies the
 boundary correspondence $\phi(J_j)=J_j'$.
  For $x_1,x_2\in E$ such that $x_1<x_2$ set
$$
 \rho_E(x_1,x_2)=\rho_E(x_2,x_1):=\mb{diam }\phi([x_1,x_2]),
 $$
 where
 $$\mb{diam }B:=\sup_{z,\z\in B}|z-\z|\quad (B\subset\C).$$
 In spite of its definition via the conformal mapping, the
 behavior  of $\rho_E$ can be characterized in purely geometrical terms as
 follows.
 According to  (\ref{1.2})
 \beq\label{1.11}
 d(J_j,E^*\setminus J_j)\ge C_5|J_j|,\quad C_5=C_2^{-1/2}.
 \eeq
 Let constant $C$ be fixed such that $0<C<\min(1,C_5/2)$. For any component
 $J_j$ of $E^*$, denote by $\tilde{J_j}$ the open interval with the
 same center of the length $(1+C)|J_j|$. For $x_1,x_2\in E$ such that
 $x_1<x_2$ consider the function
   $$\tau_{E}(x_1,x_2)=\tau_{E}(x_2,x_1)=\tau_{E,C}(x_1,x_2)$$
   $$:= \left\{\begin{array}{ll}
   \displaystyle
 \left(\frac{|J_j|}{d([x_1,x_2],J_j)}\right)^{1/2}(x_2-x_1)
 ,&\mb{ if }x_1,x_2\in\tilde{J_j}\mb{ for some }j
 \mb{ and }\\
 & x_2-x_1<d([x_1,x_2],J_j),\\[2ex]
 |J_j|^{1/2}(x_2-x_1)^{1/2},&\mb{ if }x_1,x_2\in\tilde{J_j}
 \mb{ for some }j
 \mb{ and }\\
 &d([x_1,x_2],J_j)\le x_2-x_1\le\frac{C}{2}|J_j|,\\[2ex]
 x_2-x_1,&\mb{ otherwise.}
\end{array}\right.
 $$
 \begin{th}\label{th1}
 For $x_1,x_2\in E$,
 \beq\label{1.04}
 \frac{1}{C_6}\tau_{E}(x_1,x_2)\le
 \rho_{E}(x_1,x_2)\le C_6 \tau_{E}(x_1,x_2),
 \eeq
 where $C_6=C_6(E,C)>1$.
 \end{th}
 Notice that according to Theorem \ref{th1}
 \beq\label{1.31}
 \rho_{E}(x_1,x_2)\ge C_7|x_2-x_1|\quad (x_1,x_2\in E).
 \eeq
 The main result of this paper is the following analogue of
 (\ref{1.01})-(\ref{1.02}): for $f\in BC(E)$ and $0<\al<1$,
  \beq\label{1.011}
 A_\si(f,E)=O(\si^{-\al})\quad\mb{as }\si\to\infty
 \eeq
 if and only if
\beq\label{1.021} \om_{f,E}^*(\de)=O(\de^\al)\quad \mb{as }\de\to
 +0,
 \eeq
 where
 $$\om_{f,E}^*(\de):=\sup_{x_1,x_2\in E\atop\rho_E(x_1,x_2)\le\de}|f(x_2
 )-f(x_1)|\quad (\de>0).
 $$
 The  statement (\ref{1.011})-(\ref{1.021}) follows
 immediately from the direct Theorem \ref{th2} and the inverse Theorem \ref{th3}
 below.

 Let $\om(\de),\de>0$ be a function of modulus of continuity type,
 i.e., a positive nondecreasing function with $\om(+0)=0$ such that
 \beq\label{1.5}
 \om(t\de)\le 2t\om(\de)\quad (\de>0,t>1).
 \eeq
 Denote by $BC^*_\om(E)$ the class of functions $f\in BC(E)$
 satisfying
 $$ \om_{f,E}^*(\de)\le \om(\de) \quad (\de>0).
 $$
 \begin{th}\label{th2} For $f\in BC^*_\om(E)$ and $\si\ge 1$,
 \beq\label{1.6}
 A_\si(f,E)\le
 C_8\left(\frac{||f||_{C(E)}}{\si}+\om \left(\frac{1}{\si}\right)\right).
 \eeq
 \end{th}
 \begin{th}\label{th3} Let $f\in BC(E)$ and let
 $$
 A_\si(f,E)\le\om\left(\frac{1}{\si}\right)\quad (\si\ge1).
 $$
 Then for $x_1,x_2\in
 E$,
 \beq\label{1.8}
 |f(x_2)-f(x_1)|\le C_9\, \Om(\rho_E(x_1,x_2)),
 \eeq
 where
 $$\Om(\de):=\de\left(||f||_{C(E)}+\int_\de^1\frac{\om(t)}{t^2}dt\right)
 \quad \left(0<\de\le \frac{1}{2}\right)
 $$
 and $\Om(\de):=\Om(1/2)$ for $\de>1/2$.
 \end{th}
 Let us introduce the notation we will be using throughout the
 paper.
 We continue to use the convention that $C,C_1,\ldots$ denote
 positive constants, different in different sections and depending
 only on inessential quantities. For $a,b\ge 0$ we write $a\ole b$
 if $a\le Cb$. We also write $a\asymp b$ if $a\ole b$ and $b\ole a$
 simultaneously.

 Let for $z\in\C$ and $\de>0$,
 $$D(z,\de):=\{\z:\, |\z-z|<\de\},\quad
 C(z,\de):=\{\z:\, |\z-z|=\de\},$$
 $$D(\de):=D(0,\de),\quad C(\de):=C(0,\de),
 $$
 $$\He_+:=\He ,\quad
 \He_-:=\C\setminus\ov{\He_+}.$$
 $$\N:=\{1,2,\ldots\},\quad \N_0:=\{0,1,\ldots\}.$$
 The rest of the paper is organized as follows. Since a
 significant number of proofs in this paper depends on the
 techniques for estimation of the module of path families, Section
 2 contains a brief summary of the appropriate results from
 geometric function theory. In Section 3, we compile certain facts
 about the Levin conformal mapping. In particular, the proof of
 Theorem \ref{th1} is given in this section.  Sections 4-6 present
 preliminary results for the proof of Theorem \ref{th2} given in Section
 7. Specifically, in Section 4, we construct auxiliary domains and study their
 conformal mappings onto a half-plane. Section 5 summarizes the
 relevant material on the continuous extension of functions from
 a closed subset of $\R$ into $\C$. Section 6
 provides an
 exposition of some facts from the theory of entire functions. In
 Section 8, we prove Theorem \ref{th3}.

 \absatz{Auxiliary Results about Conformal Mappings}

 In this section we discuss
 mostly known results which  concern the distortion properties
 of some conformal mappings. The results are stated in the form
 convenient for further exposition.

As usual, a Jordan curve is a continuous image of a closed
interval without intersections (except possibly endpoints). By a
curve we understand a locally rectifiable Jordan curve without
endpoints. We define a path to be the union of finitely many
mutually disjoint curves. We use $\Ga,\Ga_1,\ldots$ to denote path
families. We may use the same symbol for different families if it
does not lead to confusion. For a path family $\Ga$ denote by
$m(\Ga)$ its module, see \cite{ahl}, \cite{lehvir},  and
\cite{garmar}.   In the sequel we refer to the basic properties of
the module, such as conformal invariance, comparison principle,
composition laws, etc.
 (see monographs cited above for more details). As a rule, we
will use these properties  without further citations.

Special families of separating paths play a  useful role. Let
$G\subset\C$ be a domain. We say that a path $\ga\in G$ separates
sets $A\subset\ov{G}$ and $B\subset\ov{G}$ if $\ga$ consists of a
finite number of crosscuts of $G$ and any curve $J\subset G$
joining $A$ with $B$ has nonempty intersection with $\ga$. We
denote by $\Ga(A;B;G)$ the set of all such paths.
 Sometimes, more sophisticated families of separating curves are used.
 Let  $\tilde{G}$ be a
 compactification of $G$ by prime ends in the Carath\'eodory sense
 (see \cite{pom75}). A point $z\in G$  can also be understood as a prime
 end $Z\in\tilde{G}$ defined by a chain of concentric circles
 converging to this point.
 We say
 that a crosscut $\ga\subset G$ separates
 $Z_1,Z_2,\ldots \in\tilde{G}$ from  $\Z_1,\Z_2,
 \ldots \in\tilde{G}$ in $G$ if $G\setminus \ga$ consists of two
 connected components such that one of them is adjacent to
 $Z_1,Z_2,\ldots $ and the other to $\Z_1,\Z_2,\ldots$ (the
 adjacency    means that in the domain and the
 subdomain a prime end can be defined by the same chain of crosscuts or
  concentric circles).
 We denote by $\Ga(Z_1,Z_2,\ldots ;
 \Z_1,\Z_2,\ldots;G)$ the set of all such crosscuts.

 The  examples below state  some  well-known facts concerning
special path families.

\noindent {\bf Example 3.} For $0<\al\le 2\pi$ and $0<r_1<r_2,$
let
$$\Ga_1=\{\ga_r=\{re^{i\theta}:0<\theta<\al\pi\}:r_1<r<r_2\}.$$
 Then
 \beq\label{2.8}
 m(\Ga_1)=\frac{1}{\al\pi}\log\frac{r_2}{r_1}
 \eeq
 (see \cite[p. 8]{andbeldzj}).

 \noindent {\bf Example 4.} Let $z_1,z_2,z_3\in\C$ be distinct points
 satisfying $|z_1-z_2|\le|z_1-z_3|$. Then, for the module of
 $\Ga_2=\Ga(z_1,z_2;z_3,\infty;\C)$ we have
 \beq\label{2.3}
 \frac{1}{2\pi}\log\left|\frac{z_1-z_3}{z_1-z_2}\right|\le
 m(\Ga_2)\le
 \frac{1}{2\pi}\log\left|\frac{z_1-z_3}{z_1-z_2}\right| +2
 \eeq
 (see \cite[pp.
 98-99]{andbeldzj}).

 \noindent {\bf Example 5.} Let $z_1\in\R$ and $z_2,z_3\in \ov{\He}$
  be distinct points and let $\Ga_3=\Ga(z_1,z_2;z_3,\infty;\He)$.
  If $|z_1-z_2|\le|z_1-z_3|$ then
  \beq\label{2.6}
 \frac{1}{\pi}\log\left|\frac{z_1-z_3}{z_1-z_2}\right|\le
 m(\Ga_3)\le
 \frac{1}{\pi}\log\left|\frac{z_1-z_3}{z_1-z_2}\right| +2.
 \eeq
 If $|z_1-z_2|\ge|z_1-z_3|$ then
  \beq\label{2.7}
 m(\Ga_3)\le 2
 \eeq
 (cf.
 \cite[p. 35]{andbeldzj}).

 \noindent {\bf Example 6.} For $z\in\He$ and $\z\in\ov{\He}$ let
 $\Ga_4=\Ga(z;\z,\infty;\He).$
 Reasoning similar to that in the proof of (\ref{2.6}) and (\ref{2.7})
  demonstrates that if $\Im z<|\z-\Re
 z|$ then
  \beq\label{2.1}
 \frac{1}{\pi}\log\frac{|\z-\Re z|}{\Im z}\le
 m(\Ga_4)\le
 \frac{1}{\pi}\log\frac{|\z-\Re z|}{\Im z} +2,
 \eeq
 and if $\Re z=\z\in\R$, then
 \beq\label{2.5dop}
 m(\Ga_4)\le 2.
 \eeq
 {\bf Example 7.} Let $0<r<R$ and let curves $l_r\subset
 D(r)\cap\He$ and $l_R\subset \He\setminus D(R)$ be such that
 $$0\in l_r,\infty\in l_R,\quad C(r)\cap\ov{l_r}\neq\emptyset,
 C(R)\cap\ov{l_R}\neq\emptyset.$$
 Denote by $\Ga_5$  the family of all crosscuts of
 $G=\He\setminus(l_r\cup l_R)$ joining $l_r$ to $l_R$. Then
 \beq\label{2.16}
 m(\Ga_5)\ge\frac{\pi}{2\pi+\log\frac{R}{r}}.
 \eeq
 Indeed, let $z_r\in C(r)\cap\ov{l_r}$ and $z_R\in C(R)\cap\ov{l_R}$.
 According to (\ref{2.6}), for the family $\Ga_6$ of all crosscuts of $G$
 which separate $l_r$ and $l_R$  we have
 $$
 m(\Ga_6)\le m(\Ga(0,z_r;z_R,\infty,\He))\le
 \frac{1}{\pi}\log\frac{R}{r}+2.
 $$
 Since $m(\Ga_5)m(\Ga_6)=1$, we obtain (\ref{2.16}).

 \noindent {\bf Example 8.} Let $0<r<R$ and let $S_\pm$ be open sets
 consisting of disjoint intervals $(c_j^\pm,d_j^\pm)$ such that
 $$S_+=\cup_j(c_j^+,d_j^+)\subset(r,R),\quad
 S_-=\cup_j(c_j^-,d_j^-)\subset(-R,-r),
 $$
  \beq\label{2.27}
 d((c_j^\pm,d_j^\pm),\{0\}\cup S_\pm\setminus (c_j^\pm,d_j^\pm))\ge C_1 (d_j^\pm-
 c_j^\pm),
  \eeq
  \beq\label{2.28}
  \sum_j\left(\frac{d_j^\pm-
 c_j^\pm}{\min\{c_j^\pm,d_j^\pm\} }\right)^2\le C_2.
 \eeq
 Denote by $\Ga_7$ the family of all paths in
 $$Q=\{z\in\He:\, r<|z|<R\}
 $$
 which separate $(r,R)\setminus S_+$ from $(-R,-r)\setminus S_-$.
 \begin{lem}\label{lem2.1}
 Under the above assumptions,
 \beq\label{2.17}
 m(\Ga_7)\le\left(\pi\log\frac{R}{r}+C_3\right) \left(\log\frac{R}{r}\right)^{-2},
 \eeq
 where $C_3=C_3(C_1,C_2)$.
 \end{lem}
 {\bf Proof.}
 We  essentially follow the outline of the proof of \cite[Lemma
 4]{and08}. Let $U_j^\pm$ be the open interval with the
 center in $(c_j^\pm+d_j^\pm)/2$ and the length
 $(1+C_1/3)(d_j^\pm-c_j^\pm)$. Let
 $$V^+_j=(e_j^+,f_j^+):=(r,R)\cap U_j^+,\quad
  V^-_j=(e_j^-,f_j^-):=(-R,-r)\cap U_j^-,
  $$
  and let
  \begin{eqnarray*}
  R_j^+&:=&\left\{z=re^{i\theta}:\, e^+_j\le r\le f_j^+,0\le\theta\le
  C_4\log\frac{f_j^+}{e_j^+}\right\},\\
  R_j^-&:=&\left\{z=re^{i\theta}:\, -f^-_j\le r\le -e_j^-,
  \pi-C_4\log\frac{e_j^-}{f_j^-}\le\theta\le\pi\right\},
  \end{eqnarray*}
  where, by virtue of (\ref{2.27}), $C_4$ can be  chosen such that
  $$
  \log\frac{f_j^+}{e_j^+}\le\frac{\pi}{2C_4},\quad
  \log\frac{e_j^-}{f_j^-}\le\frac{\pi}{2C_4}.
 $$
  Consider the metrics
   $$\rho^*(z)=\left\{\begin{array}{ll}
|z|^{-1},&\mb{ if }z\in Q,\\[2ex]
 0,&\mb{ elsewhere in $\C$},
\end{array}\right.$$
  $$\rho^\pm_j(z)=\left\{\begin{array}{ll}
 C_5(f_j^\pm-e_j^\pm)^{-1},&\mb{ if }z\in R_j^\pm,\\[2ex]
 0,&\mb{ elsewhere in $\C$},
\end{array}\right.$$
 where $C_5$ is chosen such that for any $\ga\in\Ga_7$ with
 $\ov{\ga}\cap [c_j^\pm,d_j^\pm]\neq\emptyset$ we have
 $$
 \int_\ga \rho_j^\pm(z)|dz|\ge 1.
 $$
 According to our construction,
 for an arbitrary $\ga\in\Ga_7$ and
 $$
 \rho(z)=\max\left\{\rho^*(z),\sum_\pm\sum_j\rho_j^\pm(z)\left|
 \log\frac{f_j^\pm}{e_j^\pm}\right|\right\}
 $$
 we obtain
 \beq\label{2.29}
 \int_\ga \rho(z)|dz|\ge \log\frac{R}{r}.
 \eeq
 Since by  (\ref{2.28})
 \begin{eqnarray}
 &&\int_\C\rho(z)^2dm(z)\nonumber\\
 &\le& \int_\C\rho^{*}(z)^2dm(z)+
 \sum_\pm\sum_j\left(\int_\C\rho_j^{\pm}(z)^2dm(z)\right)\left|
 \log\frac{f_j^\pm}{e_j^\pm}\right|^2 \nonumber\\
 \label{2.30}
 &\le&\pi\log\frac{R}{r}+C_6,
 \end{eqnarray}
 where $dm(z)$ stands for the $2$-dimensional Lebesgue measure
 (area) on $\C$, the definition of the module (see \cite[p. 132]{lehvir}) as well as
 (\ref{2.29}) and (\ref{2.30}) yield (\ref{2.17}) with $C_3=C_6$.

 \hfill$\Box$

 Next, we cite a result by Jenkins and Oikawa \cite{jenoik} concerning Ahlfors' fundamental
inequalities.
\begin{lem}
\label{lem2.2} (\cite[inequalities (1) and (3)]{jenoik}) For
$0<r_1<r_2<\infty$, let
$$
Q=Q(r_1,r_2):=\{ re^{i\theta}:\, r_1<r<r_2,\,
-\theta_1(r)<\theta<\theta_2(r)\},
$$
where the
 functions $\theta_j,\, j=1,2$ have
  finite total variation $V_j$ on $[r_1,r_2]$
and satisfy
$$
0<\theta_0\le\theta_j(r)\le 2\pi.
$$
Then, for the module of $Q$, i.e., the module of the family
$\Ga=\Ga(Q)$ of curves separating in $Q$ its boundary circular
components, we have
 $$
\int\limits_{r_1}^{r_2}\frac{dr}{(\theta_1(r)+\theta_2(r))\, r}
\le m(\Ga)\le
\int\limits_{r_1}^{r_2}\frac{dr}{(\theta_1(r)+\theta_2(r))\, r}
+\frac{\pi}{\theta_0^2}(V_1+V_2).
 $$
\end{lem}
 In the proof of Theorem \ref{th2} we use a statement whose analogue for the case
of a domain including $\infty$ is due to Belyi (see \cite[pp.
65-66]{andbeldzj}).

Let $G\subset\C$ be a simply connected domain such that $i\in G$
and $ \infty\in L:=\partial G$.  Denote by $\Phi:G\to\He$ a
conformal mapping  satisfying $\Phi(i)=i,\Phi(\infty)=\infty$ and
let $\Psi:=\Phi^{-1}$. For $\xi\in\D^*:=\{\xi:\, |\xi|>1\}$ denote
by $\Ga_8(\xi)$ the family of all crosscuts of $\D^*$ which
separate $\xi$ from $\ov{\xi}$. The module of $\Ga_8(\xi)$
satisfies
 \beq\label{3.1b}
 m(\Ga_8(\xi))\le\frac{2}{\pi}\log\frac{4d(\xi,\{1,-1\})}{|\xi|-1}
 \eeq
 (see \cite[p. 113]{and03}).
 For $\tau\in\He,\Re\tau\neq0$ denote by $\Ga_9(\tau)$ the family
 of all crosscuts of $\He$ which separate $\tau$ from
 $-\ov{\tau}$. Conformal invariance of the module implies that
 \beq\label{3.2b}
 m(\Ga_9(\tau))=m(\Ga_8(\xi)),\quad \xi=\frac{\tau+i}{\tau-i}.
 \eeq
 Let $a=\Re \tau\neq0$ and $0<b=\Im \tau\le1/2$. Elementary
 computation involving (\ref{3.1b}) and (\ref{3.2b}) shows that
 \beq\label{3.3b}
 m(\Ga_9(\tau))\le\frac{2}{\pi}\log\frac{8\sqrt{2}(|a|+b)}{b}.
 \eeq
 \begin{lem}\label{lem3.1b}
 Let $\tau_1$ and $\tau_2$ be such that
 $0<\Im\tau_1=\Im\tau_2\le1/2,\, \Re(\tau_2-\tau_1)=2a\neq0$ and  let
 $\z_j=\Psi(\tau_j),j=1,2$. Then
 \beq\label{3.4b}
 \frac{|\z_1-\z_2|}{d(\z_1,L)}\le
 C\left(1+\frac{|a|}{b}\right)^4,\quad C=2^{14}.
 \eeq
 \end{lem}
 {\bf Proof.}
 According to (\ref{3.3b}) for $\Ga_{10}=\Ga(\tau_1;\tau_2;\He)$
 \beq\label{3.5b}
 m(\Ga_{10})\le\frac{2}{\pi}\log\frac{8\sqrt{2}(|a|+b)}{b}.
 \eeq
 In order to prove (\ref{3.4b}) we can assume that
 $$|\z_1-\z_2|>d:=d(\z_1,L).
 $$
 By virtue of (\ref{2.8}), for
 $$\Ga_{11}=\{C(\z_1,r):\, d<r<|\z_1-\z_2|\}$$
 we obtain
 \beq\label{3.6b}
 m(\Ga_{11})=\frac{1}{2\pi}\log\frac{|\z_1-\z_2|}{d}.
 \eeq
 Since any $\ga\in\Ga_{11}$ includes a subarc which belongs to
 $\Psi(\Ga_{10})$,  we have
 $$
 m(\Ga_{11})\le  m(\Psi(\Ga_{10}))=m(\Ga_{10}).
 $$
 Comparing the above inequality with (\ref{3.5b}) and (\ref{3.6b})
 we obtain (\ref{3.4b}).

 \hfill$\Box$

  \begin{lem}\label{lembel} Let
  $\Im \tau_1=1/\si,\si\ge 2,\tau_2=\tau_1+t,t\in\R,
 \z_j=\Psi(\tau_j),j=1,2$. Then for $z\in\C\setminus G$,
$$
 \left|\frac{\z_2-z}{\z_1-z}\right|\le C(1+\si|t|)^4\quad
 C=2^{14}+1.
 $$
 \end{lem}
 {\bf Proof.} Letting   $a=t/2$ and $b=1/\si$ in  Lemma \ref{lem3.1b} we
 have
 $$
  \left|\frac{\z_2-z}{\z_1-z}\right|\le \left|\frac{\z_2-\z_1}{\z_1-z}\right|
  +1\le (2^{14}+1)(1+\si|t|)^4.
  $$

 \hfill$\Box$

 The domains which appear in this paper do not
 have inner cusps on the boundary. That is,
 following \cite{tam} we say that $\C\setminus G=:K\in H$
 if any points $z,\z\in K$ can be joined by a curve
 $\ga(z,\z)\subset K$ such that
 $$|\ga(z,\z)|\le C_7|z-\z|,\quad C_7=C_7(K)\ge 1.$$
 In the case where $\infty\in G$ and $K\in H$ the geometric
 properties of
 $G$  are studied in \cite{and86}, \cite{and89}.
 In this paper we assume that $\infty\in L=\partial G$ and
 formulate some obvious modifications of results and constructions
 from these papers.

 For the rest  of this section we assume that $K\in H$ and the constants
 in the inequalities  depend only on $K$.
   Let $\tilde{G}$ be
 the compactification of $G$ by prime ends (see \cite{pom75})
 and let $\tilde{L}:=\tilde{G}\setminus G$.
 If $L$ is a Jordan curve, then $\tilde{L}=L$. Since $K\in H$ all
 $Z\in \tilde{L}$ are of first kind, i.e., they have singleton
 impressions $|Z|=z\in L$.
 For $Z\in\tilde{L}$ and $r>0$ denote by  $\ga_Z(r)=\ga_Z(r,G)
 \subset G\cap C(z,r)$ a crosscut of $G$ which
 separates $Z$ from $\infty$. For our purposes it is sufficient to
 assume that
   for any $Z\in\tilde{L}$ and $r>0$
 the crosscut $\ga_Z(r)$ is defined uniquely. We use the same symbol
 $\Phi$ to denote the homeomorphism between $\tilde{G}$ and
 $\ov{\He}$ which coincides in $G$ with the mapping $\Phi$ and let
 $\Psi=\Phi^{-1}$.

 First, note that for $Z\in\tilde{L}$ and $r>0$,
 \beq\label{2.14}
 \sup_{\z\in\ga_Z(r)}|\Phi(\z)-\Phi(Z)|\asymp \inf_{\z\in\ga_Z(r)}|\Phi(\z)-\Phi(Z)|
 \eeq
 (cf. \cite[Lemma 2]{and89}).

 If $0<r<R$ then $\ga_Z(r)$ and $\ga_Z(R)$ are the
sides of some quadrilateral $Q_Z(r,R)=Q_Z(r,R,G)\subset G$ whose
other two sides are parts of the boundary $L$. Let
$m_Z(r,R)=m_Z(r,R,G)$ be the module of this quadrilateral, i.e.,
the module of the family of curves that separate the sides
$\ga_Z(r)$ and $\ga_Z(R)$ in $Q_Z(r,R)$.
 \begin{lem}
\label{lem2.33} (\cite[Theorem 2]{and89}) Let  $Z\in \tilde{L}$
and $0<r_1<r_2<r_3$. Then
$$
0\le m_Z(r_1,r_3)-(m_Z(r_1,r_2)+m_Z(r_2,r_3))\ole 1,
$$
$$
\frac{1}{2\pi}\log\frac{r_2}{r_1}\le m_Z(r_1,r_2)\ole
 \log\frac{r_2}{r_1}+1.
$$
Moreover, for any $\z_j\in\ga_Z(r_j),\, j=1,2,$
 $$
m_Z(r_1,r_2)\le m(\Ga(Z,\z_1;\z_2,\infty;G))\le
m_Z(r_1,r_2)+C_{8}.
$$
\end{lem}
 \begin{lem}\label{lem2.3}(\cite[Lemmas 2 and 3]{and86}) Let
 $Z\in\tilde{L},\z_1,\z_2\in G,w=\Phi(Z),\tau_k=\Phi(\z_k),k=1,2$
 be such that $\Re w=\Re\tau_1=\Re \tau_2,\Im\tau_1<\Im\tau_2$. Then
 \beq\label{2.n4}
 |z-\z_1|\ole |z-\z_2|,
 \eeq
 \beq\label{2.n3}
 d(\z_k,L)\asymp|\z_k-z|,
 \eeq
 \beq\label{2.n1}
 \left|\frac{w-\tau_1}{w-\tau_2}\right|^2\ole
 \left|\frac{z-\z_1}{z-\z_2}\right|\ole
 \left|\frac{w-\tau_1}{w-\tau_2}\right|^C.
 \eeq
 \end{lem}
 The analysis of the proof of the left-hand side of (\ref{2.n1}),
 given in \cite{and86}, shows that in the case $|Z|=z\in\R$ and
 $G\subset\He$,
 it can be replaced by a sharper one
 \beq\label{2.n2}
 \left|\frac{w-\tau_1}{w-\tau_2}\right|\ole
 \left|\frac{z-\z_1}{z-\z_2}\right|.
 \eeq
 For $Z\in\tilde{L},\z\in G$, and $\de>0$ set
 $$z_\de:=\Psi(\Phi(Z)+i\de),\quad \z_\de:=\Psi(\Phi(\z)+i\de).$$
 \begin{lem}\label{lem2.4} (\cite[Lemma 4]{and86}) Under the
 above notation if $|\Phi(\z)-\Phi(Z)|\le C_9\de$ then
 $$\frac{1}{C_{10}}|z-z_\de|\le |\z-\z_\de|\le C_{10}|z-z_\de|,$$
 where $C_{10}=C_{10}(G,C_9)\ge 1$.
 \end{lem}
 For $Z\in\tilde{L},{\cal Z}\in \tilde{G}$ and $z=|Z|,\z=|{\cal Z}|$
  we denote by $r_Z({\cal Z})$ the supremum of those $r>0$
for which the arc $\ga_Z(r)$ separates $Z$ from ${\cal Z}$.
 By the definition of the class $H$ we have
 \beq\label{2.n5}
 r_Z({\cal Z})\ge C_{11}|z-\z|.
 \eeq
Moreover, if $\z\in G$ satisfies $\Re\Phi(\z)=\Re\Phi(Z)$, then
  \beq\label{2.n6}
 r_Z(\z)\le C_{12}|z-\z|
 \eeq
 (cf. \cite[p. 61]{and89}).
\begin{lem}\label{lem2.5}
 Let $x\in\R,\de>0,z\in\ov{\He}\cap C(x,\de)$, and let
 $W=\Psi(x),w=|W|,T=\Psi(z),t=|T|, w_\de=\Psi(x+i\de). $ Then
 \beq\label{2.21}
 |w-t|\ole|w-w_\de|.
 \eeq
 \end{lem}
 {\bf Proof.} If $|w-t|\ge\frac{C_{12}}{C_{11}}|w-w_\de|$ we consider
 $$\Ga_1=\Ga(x,z;x+i\de,\infty;\He),$$
 $$\Ga_2=\{C(w,r):\, C_{12}|w-w_\de|<r<C_{11}|w-t|\}.$$
 According to  (\ref{2.8}),  (\ref{2.7}), (\ref{2.n5}), and (\ref{2.n6})
 $$2\ge
 m(\Ga_1)=m(\Psi(\Ga_1))\ge m(\Ga_2)=\frac{1}{2\pi}\log
 \frac{C_{11}|w-t|}{C_{12}|w-w_\de|},
 $$
 which implies (\ref{2.21}).

 \hfill$\Box$

\absatz{The Levin Conformal Mapping}

In this section we discuss metric properties of the Levin
conformal mapping $\phi$ introduced in Section 1. We extend $\phi$
continuously to $\R$. This extension can be understood in two
different ways: either (i) for $x\in\R$ the value of $\phi(x)$ is
a point from $\partial\He_E$ or (ii)  values of $\phi$ on $\R$ are
prime ends in $\tilde{\He}_E\setminus\He_E$. Below we denote
either extension by the same letter $\Phi$.

Since by (\ref{1.1}) and (\ref{1.11}) $E$ is relatively dense with
respect to the linear Lebesgue measure, i.e., there exist $C_1$
and $C_2$ such that
 $$|[x,x+C_1]\cap E|\ge C_2\quad (x\in\R),$$
 according to \cite[Theorem 3.9]{lev93} for any component $J_j=(a_j,b_j)$ of $E^*$
 and its image $J'_j:=\phi(J_j)$ we have
 $$
 |J_j'|\ole 1
 $$
  as well as
  \beq\label{3.2}
  |\phi(z)|+1\asymp|z|+1\quad (z\in\He).
  \eeq
  The set $E$ consists of disjoint components which are either
  points or closed intervals.
  Denote by $T_j=[c_j,a_j]$ and $R_j=[b_j,d_j]$ the components of
  $E$ that  adjoin  $J_j$ and let $T_j':=\phi(T_j),
  R_j':=\phi(R_j).$
  \begin{lem}\label{lem3.1}
  The following holds:
  \beq\label{3.3}
  |J_j'|\ole|T_j'|,
  \eeq
  \beq\label{3.4}
  |J_j'|\ole|R_j'|.
  \eeq
  \end{lem}
  {\bf Proof.}
 In order to prove (\ref{3.3}) we can assume that $|J_j'|>|T_j'|$.
 Consider families of crosscuts
 $$\Ga_1=\Ga(c_j,a_j;b_j,\infty;\He),\quad \Ga_2=\{C(u_j,r)\cap\He :\,
 |T_j'|<r<|J_j'|\}.
 $$
 By virtue of (\ref{1.11}), (\ref{2.8}), and (\ref{2.6})
 \begin{eqnarray*}
 &&\frac{1}{\pi}\log\frac{|J_j'|}{|T_j'|}=m(\Ga_2)\le
  m
 (\phi(\Ga_1))\\
 &=& m(\Ga_1)
 \le \frac{1}{\pi}\log\frac{|T_j|+|J_j|}{|T_j|}+2
 \ole 1,
 \end{eqnarray*}
 which establishes the formula (\ref{3.3}).

 Similar argument applies to the proof of inequality (\ref{3.4}).

 \hfill$\Box$

  Lemma \ref{lem3.1} implies that $\C\setminus\He_E\in H$, i.e.,
 $\phi$ and its inverse mapping $\psi:=\phi^{-1}$ satisfy the
 properties mentioned in Section 2. In particular, let for $x\in
 E$ and $\de>0$,
 $$E_\de:=\{z\in\ov{\He}:\, \Im\phi (z)=\de\},\quad \rho_\de(x)
 :=d(x,E_\de).$$
 Denote by $x_\de^*\in E_\de$ any point satisfying $\Re
 x_\de^*=\Re x$. We claim that
 \beq\label{3.19}
 |x-x_\de^*|\ole\rho_\de(x).\eeq
 Indeed, let  $z\in
 E_\de$ be such that $|x-z|=\rho_\de(x)$. If
 $|x-x_\de^*|>\rho_\de(x)$, consider families of curves
 $$\Ga_1=\Ga(x,z;x_\de^*,\infty;\He),\quad \Ga_2=\{C(x,r)\cap\He :\,
 \rho_\de(x)<r<|x-x_\de^*|\}\}.
 $$
 According to (\ref{2.8}), (\ref{2.n3}), and Lemma \ref{lem2.33}
 $$\frac{1}{\pi}\log\frac{|x-x_\de^*|}{\rho_\de(x)}=m(\Ga_2)
 \le m(\Ga_1)=m(\phi(\Ga_1))\ole 1,$$
  which yields (\ref{3.19}).

 Comparing (\ref{3.19}) and  (\ref{2.n3})-(\ref{2.n2}) for $0<\de<\De$ we obtain
 \beq\label{3.20}
 \left(\frac{\de}{\De}\right)^{1/C}\ole\frac{\rho_\de(x)}{\rho_\De(x)}\asymp
 \left|\frac{x-x_\de^*}{x-x_\De^*}\right|\ole\frac{\de}{\De},
 \eeq
 where $C$ is the constant from (\ref{2.n1}).
 \begin{lem}\label{lem3.2}
 Let $x_1,x_2\in E,x_1<x_2, W_k=\phi(x_k)\in\tilde{\He}_E\setminus\He_E,
 |W_k|=w_k\in\R,k=1,2.$ Then
 \beq\label{3.6}
 \mb{\em diam }\phi([x_1,x_2])\asymp|w_k-\phi(x_k+i(x_2-x_1))|\quad
 (k=1,2).
 \eeq
 \end{lem}
 {\bf Proof.} We  only consider  the case
 $k=1$. The proof for the other case is similar.
 Let $\de=x_2-x_1,\tau_1=\phi(x_1+i\de)$. Consider the
 curve
 $$J=[x_1,x_1+i\de]\cup\left\{z=\de e^{i\theta}:\,
 0\le\theta\le\frac{\pi}{2}\right\}.
 $$
 By virtue of (\ref{2.n4}) and Lemma \ref{lem2.5}
 $$\mb{ diam }\phi(J)\ole |w_1-\tau_1|,$$
 which implies that
 \beq\label{3.61}
 B:=\mb{ diam }\phi([x_1,x_2])\ole |w_1-\tau_1|.
 \eeq
 To prove the opposite inequality, according to (\ref{2.n5}) we can
 assume that $B<r_{W_1}(\tau_1)$. Further, consider
 $$\Ga=\Ga(W_1,W_2;\tau_1,\infty;\He_E).$$
 Lemma \ref{lem2.33} and (\ref{2.7}) yield
 \begin{eqnarray*}
 \frac{1}{2\pi}\log\frac{r_{W_1}(\tau_1)}{B}&\le&
 m_{W_1}(B,r_{W_1}(\tau_1);\He_E)\\
 &\le& m(\Ga)=m(\psi(\Ga))\le2,
 \end{eqnarray*}
 which, together with (\ref{2.n5}), implies that
 \beq\label{3.62}
 B\oge r_{W_1}(\tau_1)\oge|w_1-\tau_1|.
 \eeq
 Comparing (\ref{3.61}) and (\ref{3.62}) we obtain (\ref{3.6}).

 \hfill$\Box$

 \begin{lem}\label{lem3.3}
 The following holds:
 \beq\label{3.7}
 |J_j'|\asymp|J_j|.
 \eeq
 \end{lem}
 {\bf Proof.} For sufficiently large $R\in E$ consider
 $$\Ga_1=\Ga(a_j,b_j;R,\infty;\He),\quad
 \Ga_2=\{C(u_j,r)\cap\He:\, |J_j'|<r<\z_R-u_j\},
 $$
 where ${\cal{Z}}_R:=\phi(R)$ and $\z_R:=|{\cal{Z}}_R|$.

 By (\ref{2.8}) and (\ref{2.6})
 \begin{eqnarray*}
 \frac{1}{\pi}\log\frac{\z_R-u_j}{|J_j'|}&=&m(\Ga_2)\le
 m(\phi(\Ga_1))\\
 &=&m(\Ga_1)\le\frac{1}{\pi}\log\frac{R-a_j}{|J_j|}+2,
 \end{eqnarray*}
 which implies that
 $$
 \frac{|J_j|}{|J_j'|}\ole\frac{R-a_j}{\z_R-u_j}.
 $$
 Letting $R\to\infty$ and using (\ref{3.2}) we
 have
 \beq\label{3.8}
 |J_j'|\oge|J_j|.
 \eeq
 In order to prove the opposite inequality, for sufficiently large
 $R$ consider the family $\Ga_3$ of all crosscuts of
 $\He\setminus(J_j'\cup\{u_j+iv:\, v\ge R\})$ joining $J_j'$ with
 $\{u_j+iv:\, v\ge R\}$. Let $\Ga_4=\{\ga\cap\He_E:\,
 \ga\in\Ga_3\}.$
 According to (\ref{2.16})
 \beq\label{3.9}
 m(\Ga_4)\ge m(\Ga_3)\ge
 \pi\left(\log\frac{C_3R}{|J_j'|}\right)^{-1}.
 \eeq
 Moreover, by virtue of (\ref{3.2}) there is a constant $C_4$
  such that for sufficiently large $R$
 $$\left|\psi(u_j+iv)-\frac{b_j-a_j}{2}\right|\ge C_4R\quad
 (v\ge R).$$
 Let $\Ga_5$ be the family of all paths in
 $$Q=\left\{z\in\He:\,
 \frac{|J_j|}{2}<\left|z-\frac{b_j+a_j}{2}\right|
 <C_4R\right\}
 $$
 which separate $(r,R)\setminus E^*$ from $(-R,-r)\setminus E^*$.
 By  (\ref{1.2}) and Lemma \ref{lem2.1}
 \beq\label{3.10}
 m(\Ga_4)=m(\psi(\Ga_4))\le m(\Ga_5)\le
 \pi\left(\log\frac{C_5R}{|J_j|}\right)
 \left(\log\frac{2C_4R}{|J_j|}\right)^{-2}.
 \eeq
 Comparing (\ref{3.9}) and (\ref{3.10}) and letting
  $R\to\infty$ we obtain
 \beq\label{3.101}
 |J_j'|\ole |J_j|,
 \eeq
 which gives (\ref{3.7}) when combined with (\ref{3.8}).

 \hfill$\Box$

 For $x\in E$ and $\de>0$ let $W:=\phi(x),w:=|W|,
 w_\de:=\phi(x+i\de).$ The next three lemmas state estimates for
 the quantity $|w-w_\de|$.
 \begin{lem}\label{lem3.4}
 The following holds:
 \beq\label{3.11}
 |w-w_\de|\oge\de.
 \eeq
 \end{lem}
 {\bf Proof.}
 For sufficiently large $R\in E$ consider
 $$\Ga_1=\Ga(x,x+i\de;R,\infty,\He),$$
 $$
 \Ga_2=\{ C(w,r)\cap\He:\, C_{12}|w-w_\de|<r<C_{11}|w-\z_R|\},
 $$
 where $C_{11}$ and  $C_{12}$ are the constants from (\ref{2.n5}) and
 (\ref{2.n6}); and $\z_R$ is the impression of the prime end
 $\phi(R)$. According to (\ref{2.8}) and (\ref{2.6})
 \begin{eqnarray*}
 &&2+\frac{1}{\pi}\log\frac{R-x}{\de}\ge m(\Ga_1)= m(\phi(\Ga_1))\\
 &\ge& m(\Ga_2)\ge\frac{1}{\pi}\log\frac{C_{11}|\z_R-w|}{C_{12}
 |w_\de-w|},
 \end{eqnarray*}
 from which, after taking into account (\ref{3.2}) and letting
   $R\to\infty$, we obtain (\ref{3.11}).

  \hfill$\Box$

 Let $\tilde{J_j}$ be as in Section 1, i.e., it is the open
 interval with the same center as $J_j$ and the length
 $(1+C)|J_j|$ where $C$ is sufficiently small.
 \begin{lem}\label{lem3.5}
 For $x\in \R\setminus\cup_j\tilde{J_j}$,
 \beq\label{3.12}
 |w-w_\de|\ole\de.
 \eeq
 \end{lem}
 {\bf Proof.}
 Essentially, we have to mimic the proof of (\ref{3.101}). Hence,
 we only sketch it. For sufficiently large $R$ let
 $$S=\{\tau=\phi(x+it):\, 0\le t\le\de\},$$
 $$T=T_R=\{\tau=\phi(x+it):\, t\ge R\}.$$
 Consider the family $\Ga_1$ of curves joining $S$ and $T$ in
  $\He\setminus (S\cup T)$. Let
 $$\Ga_2=\{\ga\cap\He_E:\, \ga\in\Ga_1\}.$$
 According to
 (\ref{2.16}), (\ref{2.n4}), and (\ref{3.2})
 \beq\label{3.13}
 m(\Ga_2)\ge
 m(\Ga_1)\ge\pi\left(\log\frac{C_6R}{|w-w_\de|}\right)^{-1}.
 \eeq
 Furthermore, by virtue of  (\ref{1.2}) and Lemma \ref{lem2.1}
 for the module of the family $\Ga_3$ of all paths in
 $$Q=\{z\in\He:\, \de<|z-x|<R\}
 $$
 which separate $(\de,R)\setminus E^*$ from $(-R,-\de)\setminus
 E^*$ we have
 \beq\label{3.14}
 m(\Ga_3)\le\pi\left(\log\frac{C_{7}R}{\de}\right)
 \left(\log\frac{R}{\de}\right)^{-2}.
 \eeq
 Since
 $$m(\Ga_2)= m(\psi(\Ga_2))\le m(\Ga_3),$$
 comparing (\ref{3.13}), (\ref{3.14})  and letting
 $R\to\infty$ we obtain  (\ref{3.12}).

  \hfill$\Box$

 \begin{lem}\label{lem3.6}
 Let $x\in\tilde{J_j}$ for some $j$.

 (i) If $\de\le d(x,J_j)$, then
 \beq\label{3.15}
 |w-w_\de|\asymp\left(\frac{|J_j|}{d(x,J_j)}\right)^{1/2}\de.
 \eeq

 (ii) If $d(x,J_j)<\de\le |J_j|$, then

 \beq\label{3.161}
 |w-w_\de|\asymp|J_j|^{1/2}\de^{1/2}.
 \eeq

 (iii) If $\de\ge |J_j|$, then
 \beq\label{3.17}
 |w-w_\de|\asymp\de.
 \eeq
 \end{lem}
 {\bf Proof.}
 Without loss of generality we assume that $x\ge b_j$, i.e., $x\in
 R_j=[b_j,d_j]$. First, we provide some auxiliary estimates.
 According to the definition of $\tilde{J_j}$
 \beq\label{3.181}
 x-b_j\le C|J_j|<\frac{1}{2}|R_j|.
 \eeq
 We claim that
 \beq\label{3.18}
 w-u_j\ole v_j=|J_j'|.
 \eeq
 Indeed, by (\ref{2.6}), for the module of the family
 $$\Ga_1=\Ga(a_j,b_j;x,\infty;\He)$$
 we have
 \beq\label{3.191}
 m(\Ga_1)\le\frac{1}{\pi}\log\frac{x-a_j}{b_j-a_j}+2\ole 1.
 \eeq
 on the other hand, if $w-u_j>v_j$, by virtue of (\ref{2.8})
 \beq\label{3.201}
 m(\Ga_1) = m(\phi(\Ga_1))\ge
 m(\Ga_2)=\frac{1}{\pi}\log\frac{w-u_j}{v_j},
 \eeq
 where
 $$\Ga_2=\{C(u_j,r)\cap\He:\,v_j<r<w-u_j\}.$$
 Comparing (\ref{3.191}) and (\ref{3.201}) we obtain (\ref{3.18}).

 Furthermore,
 \beq\label{3.21}
 w-u_j\asymp(x-b_j)^{1/2}|J_j|^{1/2}.
 \eeq
 Indeed, by (\ref{2.6}) and (\ref{3.181}), for the module of the family
 $$\Ga_3=\Ga(x,b_j;a_j,\infty;\He)$$
 we have
 \beq\label{3.22}
 \frac{1}{\pi}\log\frac{|J_j|}{x-b_j}\le m(\Ga_3)\le
 \frac{1}{\pi}\log\frac{|J_j|}{x-b_j} +C_{8}.
 \eeq
 On the other hand, making use of Lemmas \ref{lem2.2} and
 \ref{lem2.33}, we have
 \beq\label{3.23}
 \frac{2}{\pi}\log\frac{|J_j'|}{w-u_j}-C_{9}\le m(\phi(\Ga_3))\le
 \frac{2}{\pi}\log\frac{|J_j'|}{w-u_j} +C_{9},
 \eeq
  which, together with (\ref{3.7}), (\ref{3.18}), and (\ref{3.22}),
   implies (\ref{3.21}).

  We are now in a position  to prove (\ref{3.15})-(\ref{3.17}).
  We compare the modules of pares of conformally invariant
  families of curves in $\He$ and $\He_E$. Since the reasoning
  mimics the proof of (\ref{3.21}) we only indicate the appropriate
  choice of the families of crosscuts leaving the details to the
  reader.

  (i) In this case, by (\ref{2.14}) and Lemma \ref{lem2.3},
  $|w-w_\de|\ole w-u_j$ and for the family
  $$\Ga_4=\Ga(x,x+i\de;b_j,\infty;\He)$$
  according to (\ref{2.6}) and (\ref{2.7}) we have
  $$\frac{1}{\pi}\log\frac{x-b_j}{\de}\le m(\Ga_4)\le
 \frac{1}{\pi}\log\frac{x-b_j}{\de} +2,
 $$
 as well as
 $$
 \frac{1}{\pi}\log\frac{w-u_j}{|w-w_\de|}-C_{10}\le m(\phi(\Ga_4))\le
 \frac{1}{\pi}\log\frac{w-u_j}{|w-w_\de|} +C_{10}.
 $$
 Therefore,
 $$|w-w_\de|\asymp\frac{w-u_j}{x-b_j}\de,
 $$
 which,  when combined with (\ref{3.21}), gives (\ref{3.15}).

 (ii)
 In this case, by (\ref{2.14}) and Lemma \ref{lem2.3},
 $w-u_j\ole |w-w_\de|\ole|J_j'|$ and by virtue of (\ref{2.6})
 for the family
  $$\Ga_5=\Ga(x,x+i\de;a_j,\infty;\He)$$
   we obtain
  $$\frac{1}{\pi}\log\frac{|J_j|}{\de}\le m(\Ga_5)\le
 \frac{1}{\pi}\log\frac{|J_j|}{\de} +2.
 $$
 Since, by Lemmas \ref{lem2.2} and \ref{lem2.33},
 $$
 \frac{2}{\pi}\log\frac{|J_j'|}{|w-w_\de|}-C_{11}\le m(\phi(\Ga_5))\le
 \frac{2}{\pi}\log\frac{|J_j'|}{|w-w_\de|} +C_{11},
 $$
 the two above double inequalities,
  together with (\ref{3.7}), imply (\ref{3.161}).

  (iii) The part $|w-w_\de|\oge\de$ follows from (\ref{3.11}). In
  order to prove the opposite inequality $|w-w_\de|\ole\de$ we
  have to repeat  the proof of (\ref{3.101}) word for word.

  \hfill$\Box$

  \noindent {\bf Proof of Theorem \ref{th1}.} According to Lemma
  \ref{lem3.2}
  $$\rho_E(x_1,x_2)\asymp|w_1-\tau_1|,
  $$
  where $w_1=|W_1|,W_1=\phi(x_1)$, and
  $\tau_1=\phi(x_1+i|x_2-x_1|)$. Therefore, Lemmas
   \ref{lem3.4}-\ref{lem3.6} imply (\ref{1.04}).

   \hfill$\Box$

 \absatz{Auxiliary Domains and Their Conformal Mappings}

 In this section, starting with $E$ and $E^*=\cup_jJ_j=\cup_j
 (a_j,b_j)$, we construct two auxiliary domains
 $G_\pm\supset\He_\pm$ and study their conformal mappings onto
 $\He_\pm$.

 Consider the curves
 \begin{eqnarray*}
 S_j^+&:=&[a_j-2it_j,b_j-2it_j]\cup\{z:\,|z-(a_j-it_j)|=t_j,\Re
 z\le a_j\}\\
 &&\cup \{z:\,|z-(b_j-it_j)|=t_j,\Re
 z\ge b_j\},
 \end{eqnarray*}
 where $t_j:=C_5|J_j|/3$ and $C_5$ is the constant from
 (\ref{1.11}). Denote by $G_+\supset \He_+$ the Jordan domain
 bounded by
 $$L^+=\partial G_+=E\cup(\cup_jS_j^+)$$
 and let
 $$G_-:=\{z:\, \ov{z}\in G_+\},\quad L^-:=\partial G_-,
 \quad S_j^-:=\{z:\, \ov{z}\in S_j^+\}.$$
 Denote by $\phi_\pm:\, G_\pm\to\He_\pm$ the conformal mapping
 normalized by
 $$\phi_\pm(\infty)=\infty,\quad \phi_\pm(\pm i)=\pm i,$$
 and let $\psi_\pm:=\phi_\pm^{-1}$ be the inverse mapping. By the
 symmetry
 $$\phi_+(x)=\phi_-(x)\quad (x\in E).
 $$
 Since, by (\ref{1.11})  $\C\setminus G_\pm\in H$, the conformal
 mappings $\phi_\pm$ and $\psi_\pm$ possess appropriate
 properties stated in Section 2.
 \begin{lem}\label{lem4.1} For $z\in G_\pm$,
 \beq\label{4.1}
 |\phi_\pm(z)|+1\asymp |z|+1.
 \eeq
 \end{lem}
 {\bf Proof.} Let $w=w^\pm=\phi_\pm(z)$. Without loss of generality we
 can assume that $z\in\He_\pm, |z|>2d( i,E)$, and $|w|>1$.
 According to (\ref{2.1}), for $\Ga_1=\Ga(\pm i;z,\infty;G_\pm)$ and
 $\Ga_2=\Ga(\pm i;z,\infty;\He_\pm)$ we
 have
 \beq\label{4.2}
 \frac{1}{\pi}\log|w|\le m(\phi_\pm(\Ga_1))\le
 \frac{1}{\pi}\log|w|+2,
 \eeq
 as well as
 \beq\label{4.3}
 m(\Ga_1)\le m(\Ga_2)\le \frac{1}{\pi}\log
 |z|+2.
 \eeq
 Comparing (\ref{4.3}) and the left-hand side of  (\ref{4.2}) we
 obtain
 \beq\label{4.4}
 |w|\ole|z|.
 \eeq
 Furthermore, let $x\in E$ be such that $|x\mp i|=d(i,E)=:d$
 and let
 $$\Ga_3=\{\ga_x(r,G_\pm):\, d<r<|z|-d\}.
 $$
 Notice that
 \begin{eqnarray*}
 m(\Ga_1)&\ge&
 m(\Ga_3)=\int_{d}^{|z|-d}\frac{dr}{|\ga_x(r,G_\pm)|}\\
 &\ge&\int_d^{|z|/2}\frac{dr}{\pi
 r}-\frac{1}{\pi^2}\int_d^{|z|}\frac{|\ga_x(r,G_\pm)|-\pi
 r}{r^2}dr.
 \end{eqnarray*}
 Since, by our assumption (\ref{1.2})
 $$\int_d^{\infty}\frac{|\ga_x(r,G_\pm)|-\pi
 r}{r^2}dr\ole 1,
 $$
 we have
 \beq\label{4.44}
 m(\Ga_1)\ge\log|z|-C_1,
 \eeq
 which, together with the right-hand side of (\ref{4.2}), implies
 that
 $$|w|\oge|z|.$$
 Comparing the above inequality with (\ref{4.4}) we obtain
 (\ref{4.1}).

 \hfill$\Box$

  Denote by $T_j=[c_j,a_j]$ and
 $R_j=[b_j,d_j]$ the components of $E$ adjoint to $J_j$ and let
 $$S_j^*:=\phi_\pm(S_j^\pm),\quad T_j^*:=\phi_\pm(T_j),\quad
 R_j^*:=\phi_\pm(R_j).
 $$
 \begin{lem}\label{lem4.2}
 The following holds:
 \beq\label{4.5}
 |S^*_j|\ole|T^*_j|,
 \eeq
 \beq\label{4.6}
 |S^*_j|\ole|R^*_j|,
 \eeq
 \beq\label{4.7}
 |S^*_j|\asymp|J_j|.
 \eeq
 \end{lem}
 {\bf Proof.} By virtue of (\ref{1.11}), (\ref{2.3}), and
 (\ref{2.6}) for
 $$\Ga_1=\Ga(a_j,b_j;c_j,\infty;G_\pm),\quad
 \Ga_2=\Ga(a_j,b_j;c_j,\infty;\C)
 $$
 we have
 \begin{eqnarray*}
 1&\ole& \frac{1}{2\pi}\log\frac{b_j-c_j}{b_j-a_j}\le m(\Ga_2)\le
 m(\Ga_1)\\
 &=&
 m(\phi_\pm(\Ga_1))\le\frac{1}{\pi}\log\frac{|S_j^*|+|T_j^*|}{|S_J^*|}+2,
 \end{eqnarray*}
 from which (\ref{4.5}) follows.

 The same reasoning applies to the proof of (\ref{4.6}).

 In order to prove (\ref{4.7}),  for sufficiently large
 $R$  consider the points $R\in E, t_R:=\phi_\pm(R)$ and families of curves
 $$\Ga_3=\Ga(a_j,b_j;R,\infty;G_\pm),\quad
 \Ga_4=\Ga(a_j,b_j;R,\infty;\He_\pm),
 $$
 $$\Ga_5=\{\ga_{a_j}(r,G_\pm):\, |J_j|<r<R-a_j\}.
 $$
 According to (\ref{2.6}) we obtain
 \beq\label{4.8}
 \frac{1}{\pi}\log\frac{t_R-\phi_\pm(a_j)}{|S_j^*|}\le
 m(\phi_\pm(\Ga_3))\le
 \frac{1}{\pi}\log\frac{t_R-\phi_\pm(a_j)}{|S_j^*|} +2,
 \eeq
  \beq\label{4.9}
  m(\Ga_3)\le  m(\Ga_4)\le \frac{1}{\pi}\log\frac{R-a_j}{|J_j|}+2.
  \eeq
  Moreover, using the assumption (\ref{1.2}) and reasoning as in
  the proof of (\ref{4.44}) we have
  \beq\label{4.10}
  m(\Ga_3)\ge m(\Ga_5)\ge\frac{1}{\pi}\log\frac{R-a_j}{|J_j|}-C_2.
  \eeq
  The inequalities (\ref{4.8})-(\ref{4.10}) yield
  $$\frac{|S_j^*|}{|J_j|}\asymp
  \frac{t_R-\phi_\pm(a_j)}{R-a_j}.
  $$
  Letting  $R\to \infty$
  and applying (\ref{4.1}) we obtain (\ref{4.7}).

  \hfill$\Box$

  For $x\in E,\z\in\ov{G_\pm}$, and $\si>0$ set
  $$\z_\si^\pm:=\psi_\pm\left(\phi_\pm(\z) \pm\frac{i}{\si}\right),
  $$
  $$
  L^\pm_\si:=\left\{\z\in G_\pm:\,
  \Im\phi_\pm(\z)=\pm\frac{1}{\si}\right\},\quad
  d_\si(x):=d(x,L^\pm_\si\cap\He_\pm).
  $$
  \begin{lem}\label{lem4.4}
  For $x\in E$, $\z\in\He_\pm\cap D(x,d_\si(x))$, and $\si\ge 1$ we have
  \beq\label{4.15}
  |x^\pm_\si-x|\ole d_\si(x),
  \eeq
  \beq\label{4.16}
  |\z^\pm_\si-\z|\asymp d_\si(x),
  \eeq
  \beq\label{4.14}
  |\z^\pm_\si-x|\oge d_\si(x).
  \eeq
  \end{lem}
  {\bf Proof.}
  According to (\ref{2.n5}), in order to prove (\ref{4.15}) we can
  assume that $d_\si(x)< r_x(x^\pm_\si,G_\pm)=:B$.
  Consider
   $$\Ga_1=\Ga(x,z;x_\si^\pm,\infty;G_\pm),\quad
   \Ga_2=\{ C(x,r):\, d_\si(x)<r<B\},$$
   where $z\in L_\si^\pm\cap\ov{\He_\pm}$ is such that
   $|x-z|=d_\si(x)$.
   By virtue of (\ref{2.8}) and  (\ref{2.7})
   $$\frac{1}{2\pi}\log\frac{B}{d_\si(x)}=m(\Ga_2)\le
   m(\Ga_1)=m(\phi_\pm (\Ga_1))\le 2,$$
   which, together with (\ref{2.n5}), implies  (\ref{4.15}).

   The inequality (\ref{4.16}) follows immediately from Lemma \ref{lem2.4} and
   (\ref{4.15}).

   Furthermore, by Lemma \ref{lem2.3}
   \beq\label{4.pom}
   |\z_\si^\pm-\z|\ole d(\z_\si^\pm,L^\pm)\le |\z_\si^\pm-x|,
   \eeq
   which gives (\ref{4.14}) when combined with (\ref{4.16}).

   \hfill$\Box$

  \begin{lem}\label{lem4.5}
 Let $x\in E,\z\in\He_\pm,|\z-x|\ge d_\si(x),\si\ge 1.$ Then
 \beq\label{4.18}
 \left|\frac{\z^\pm_\si-\z}{\z^\pm_\si-x}\right|
 \ole\left|\frac{d_\si(x)}{\z-x}\right|^{1/4}.
 \eeq
 \end{lem}
 {\bf Proof.} Let $t=\phi_\pm(x), \tau^\pm:=\phi_\pm(\z),
 \tau^\pm_\si:=\tau^\pm\pm i/\si$.
  Let $\Ga_1$ be the family of all closed curves in $G_\pm$
  which separate $\z$ and $\z_\si^\pm$ from  $L^\pm$.
 Consider
 $$\Ga_2=\Ga(\z,\z_\si^\pm;x,\infty;G_\pm)\cup\Ga_1,\quad
 \Ga_3=\Ga(x,x_\si^\pm;\z,\infty;G_\pm).$$
 According to (\ref{2.3}),  (\ref{2.1}), and (\ref{4.pom}) we have
 \beq\label{4.19}
 m(\Ga_2)\le\frac{1}{\pi}\log\left|\frac{\z^\pm_\si-x}{\z^\pm_\si-\z}\right|
 +C_3.
 \eeq
 Moreover, setting
 $$\Ga_4=\left\{C(\tau_\si^\pm,r):\, \frac{1}{\si}<r<|\tau_\si^\pm-t|\right\}
 $$
 by virtue of (\ref{2.8}) we obtain
 $$
 m(\Ga_2)=m(\phi_\pm(\Ga_2))\ge m(\Ga_4)=\frac{1}{2\pi}\log
 \si|\tau_\si^\pm-t|,
 $$
 which, together with (\ref{4.19}), yields
 \beq\label{4.20}
 \left|\frac{\z^\pm_\si-\z}{\z^\pm_\si-x}\right|\ole(\si|\tau_\si^\pm-t|)^{-1/2}.
 \eeq
 Next, setting
 $\Ga_5=\Ga(x,x_\si^\pm;\z,\infty;\C)$ and applying (\ref{2.3}) we have
 \beq\label{4.21}
 m(\Ga_3)\ge
 m(\Ga_5)\ge\frac{1}{2\pi}\log\left|\frac{x-\z}{x-x_\si^\pm}\right|.
 \eeq
 Furthermore, since by (\ref{2.14})
 $|\tau^\pm-t|\oge1/\si$ according to (\ref{2.6}) and (\ref{2.7})
 we obtain
 $$
 m(\phi_\pm(\Ga_3))\le\frac{1}{\pi}\log \si|\tau_\si^\pm-t| +C_4,
 $$
 which, together with (\ref{4.21}), implies that
 \beq\label{4.22}
 \left|\frac{x-x_\si^\pm}{x-\z}\right|\oge
 (\si|\tau_\si^\pm-t|)^{-2}.
 \eeq
 Comparing (\ref{4.15}), (\ref{4.20}), and  (\ref{4.22}) we obtain
 (\ref{4.18}).

 \hfill$\Box$

 \begin{lem}\label{lem4.6}
 Let $z_1,z_2\in\He_\pm$ be such that
 $$\Re z_1=\Re z_2=x\in E,\quad \Im z_2>\Im z_1,$$
 and let $t=\phi_\pm(x),w_k^\pm=\phi_\pm(z_k),k=1,2.$ Then
 \beq\label{4.231}
 |w_k^\pm-t|\ole|\Im w^\pm_k|,
 \eeq
 \beq\label{4.23}
 \left|\frac{w_2^\pm-t}{w_1^\pm-t}\right|\ole
 \left|\frac{z_2-x}{z_1-x}\right|,
 \eeq
 \beq\label{4.201}
 |w_1^\pm-t|\ole |w_2^\pm-t|.
 \eeq
 \end{lem}
 {\bf Proof.}
 Let $w^\pm_k=x_k\pm iy_k$. If $|t-x_k|>y_k$ consider
 $$\Ga_1=\Ga_{1,k}=\{C(x_k,r)\cap\He:\, y_k<r<|t-x_k|\},$$
 $$\Ga_2=\Ga_{2,k}=\Ga(z_k;x,\infty;G_\pm),\quad
 \Ga_3=\Ga_{3,k}=\Ga(z_k;x,\infty;\He_\pm).
 $$
 Since, by (\ref{2.8}) and (\ref{2.5dop}),
 $$
 \frac{1}{\pi}\log\frac{|t-x_k|}{y_k}= m(\Ga_1)
 =m(\psi_\pm(\Ga_1))\le m(\Ga_2)\le m(\Ga_3)\le 2,
 $$
 we have $|t-x_k|\ole y_k$, which establishes the formula (\ref{4.231}).

 Moreover, let
 $$\Ga_4=\Ga(x,z_1;z_2,\infty;G_\pm),\quad
 \Ga_5=\Ga(x,z_1;z_2,\infty;\He_\pm).
 $$
 According to (\ref{2.6})
 \begin{eqnarray*}
 \frac{1}{\pi}\log
 \left|\frac{w_2^\pm-t}{w_1^\pm-t}\right|& \le &
 m(\phi_\pm (\Ga_4))= m(\Ga_4)\\ &\le &m(\Ga_5)
 \le  \frac{1}{\pi}\log \left|\frac{z_2-x}{z_1-x}\right| +2,
 \end{eqnarray*}
 which implies (\ref{4.23}).

 Next, if $|t-w_1^\pm|>|t-w_2^\pm|$ we consider
 $$\Ga_6=\Ga(x,z_2;z_1,\infty;G_\pm),\quad
 \Ga_7=\Ga(x,z_2;z_1,\infty;\He_\pm),
 $$
 $$
 \Ga_8=\{C(t,r)\cap\He_\pm:\, |t-w_2^\pm|<r<|t-w_1^\pm|\}.$$
 By virtue of (\ref{2.8}) and (\ref{2.7}) we have
 \begin{eqnarray*}
 \frac{1}{\pi}\log\left|\frac{t-w_1^\pm}{t-w_2^\pm}\right|&=&
 m(\Ga_8)\le m(\phi_\pm(\Ga_6))\\
 &=&m(\Ga_6)\le m(\Ga_7)\le 2,
 \end{eqnarray*}
 from which (\ref{4.201}) follows.

  \hfill$\Box$

  \begin{lem}\label{lem4.3}
  For $z\in\He_\pm$ with $d(z,E)\le 3$,
  \beq\label{4.11}
  |\Im\phi_\pm(z)|\ole 1,
  \eeq
  \beq\label{4.12}
  |\phi'_\pm(z)|\oge 1.
  \eeq
  \end{lem}
  {\bf Proof.}
  Let $w=w^\pm=\phi_\pm(z)$ and let $R\in E$ be sufficiently
  large. Let $x\in E$ satisfy $|z-x|=d(z,E)$,
  $t_R=\phi_\pm(R)$, and let
  $$\Ga_1=\Ga(z;R,\infty;G_\pm),\quad
  \Ga_2=\{\ga_x(r,G_\pm):\, 3<r<R-x\}.
  $$
  By virtue of (\ref{2.1})
  we have
  \beq\label{4.13}
  m(\Ga_1)=m(\phi_\pm(\Ga_1))\le\frac{1}{\pi}\log\frac{t_R-\Re w}{\Im w}+2.
  \eeq
  Moreover, reasoning as in the proof of (\ref{4.44}) we obtain
  \beq\label{4.131}
  m(\Ga_1)\ge m(\Ga_2)=\int_3^{R-x}\frac{dr}{|\ga_x(r,G_\pm)|}
  \ge\frac{1}{\pi}\log(R-x)-C_5.
  \eeq
 Therefore,  (\ref{4.13}) and  (\ref{4.131}) yield
 $$|\Im w|\ole\frac{t_R-\Re w}{R-x}.$$
 Applying (\ref{4.1}) and letting $R\to\infty$ we
 have (\ref{4.11}).

 Next, let $w=u\pm iv,\xi=\xi_\pm=\psi_\pm(u)$. Since by (\ref{2.n3})
 $|z-\xi|\asymp d(z,E),
 $
 according to Lemma \ref{lem2.33} for
 $\Ga_3=\Ga(\xi,z;R,\infty;G_\pm)$ we obtain
 \begin{eqnarray*}
 m(\Ga_3)&\le& m_\xi(|\xi-z|,|\xi-R|;G_\pm)+C_6\\
 &\le&\frac{1}{\pi}\log\left|\frac{\xi-R}{\xi-z}\right|+C_6.
 \end{eqnarray*}
 In the opposite direction, by virtue of (\ref{2.6})
 $$m(\Ga_3)=m(\phi_\pm(\Ga_3))\ge\frac{1}{\pi}\log\frac{t_R-u}{v}.
 $$
 Therefore,
 $$
 \frac{v}{d(z,E)}\oge\frac{t_R-u}{|R-\xi|}.
 $$
 Letting $R\to\infty$ and taking into account
 (\ref{4.1}) we have
 $$\frac{|\Im w|}{d(z,E)}\oge 1.
 $$
 To complete the proof of (\ref{4.12}) we have to use the
 immediate consequence of the Koebe $1/4$-Lemma, i.e., the
 inequality
 $$|\phi_\pm'(z)|\ge\frac{1}{4} \frac{|\Im w|}{d(z,L^\pm)}
 \ge\frac{1}{4} \frac{|\Im w|}{d(z,E)}
 $$
 (see \cite[p. 58]{andbeldzj}).

  \hfill$\Box$

  \begin{lem}\label{lem4.7}
  Let $x\in E,0<\de\le 1,t=\phi_\pm(x),\tilde{t}_\de^\pm=\phi_\pm(x\pm
  i\de).$ Then
  \beq\label{4.24}
  |t-\tilde{t}^\pm_\de|\oge\de.
  \eeq
  \end{lem}
  {\bf Proof.}
  Let $R$ be sufficiently large and let $t_R:=\phi_\pm(R)$.
  According to
  (\ref{2.6})
  for
  $$\Ga_1=\Ga(x,x\pm i\de;R,\infty;G_\pm),\quad
  \Ga_2=\Ga(x,x\pm i\de;R,\infty;\He_\pm)
  $$
  we have
  \begin{eqnarray*}
  \frac{1}{\pi}\log\frac{t_R-t}{|\tilde{t}_\de^\pm-t|}&\le&
  m(\phi_\pm(\Ga_1))=m(\Ga_1)\le m(\Ga_2)\\
  &\le&\frac{1}{\pi}\log\frac{R-x}{\de}+2,
  \end{eqnarray*}
  which yields
  $$\frac{\de}{|\tilde{t}^\pm_\de-t|}\ole\frac{R-x}{t_R-t}.
  $$
  Letting $R\to\infty$ and applying (\ref{4.1}) we
  obtain (\ref{4.24}).

   \hfill$\Box$

 Next,  we improve the inequality (\ref{4.24})
   for points $x\in E$ close to the components $J_j=(a_j,b_j)$ of
   $E^*$.
   \begin{lem}\label{lem4.8}
   Under the assumptions and notation of Lemma \ref{lem4.7} for $x\in
   E$ such that $d(x,J_j)\le C_5|J_j|/2$, where $C_5$ is the
   constant from (\ref{1.11}), the following inequalities hold.

   (i) If $\de<d(x,J_j)$, then
   \beq\label{4.25}
   |t-\tilde{t}^\pm_\de|\oge\left(\frac{|J_j|}{d(x,J_j)}\right)^{1/2}\de.
   \eeq

   (ii) If $d(x,J_j)\le\de\le |J_j|$, then
   \beq\label{4.26}
   |t-\tilde{t}^\pm_\de|\oge |J_j|^{1/2}\de^{1/2}.
   \eeq
   \end{lem}
   {\bf Proof.}
   Let $S^*_j:=\phi_\pm(S_j^\pm)=(a_j^*,b_j^*)$. Without loss of
   generality we can assume that $x\ge b_j$.

   (i) First, we claim that
   \beq\label{4.27}
   |t-b^*_j|\oge(x-b_j)^{1/2}|J_j|^{1/2}.
   \eeq
   Indeed, by Lemmas \ref{lem2.2} and \ref{lem2.33} for the module
   of the family
   $\Ga_1=\Ga(x,b_j;a_j,\infty;G_\pm)$
   we have
   \beq\label{4.28}
   m(\Ga_1)\le m_x(x-b_j,x-a_j;G_\pm)+C_7\le
   \frac{1}{2\pi}\log\frac{|J_j|}{x-b_j}+C_{8}.
   \eeq
   Moreover, according to (\ref{2.6})
   \beq\label{4.29}
   m(\phi_\pm(\Ga_1))\ge\frac{1}{\pi}\log\frac{t-a_j^*}{t-b_j^*}
   \ge\frac{1}{\pi}\log\frac{|S_j^*|}{t-b_j^*}.
   \eeq
   Therefore, (\ref{4.7}), (\ref{4.28}), and  (\ref{4.29}) imply
   (\ref{4.27}).

   Next, by virtue of (\ref{2.6}) for
   $$  \Ga_2=\Ga(x,x\pm i\de;b_j,\infty;G_\pm),\quad
    \Ga_3=\Ga(x,b_j;a_j,\infty;\He_\pm)$$
    we have
    \begin{eqnarray*}
    \frac{1}{\pi}\log\frac{t-b_j^*}{|t-\tilde{t}^\pm_\de|}&\le& m(\phi_\pm
    (\Ga_2))=m(\Ga_2)\le m(\Ga_3)\\
    &\le& \frac{1}{\pi}\log\frac{x-b_j}{\de}+2,
    \end{eqnarray*}
    which gives  (\ref{4.25}) when combined with (\ref{4.27}).

    (ii) By Lemmas \ref{lem2.2} and \ref{lem2.33} for
    $$  \Ga_4=\Ga(x,x\pm i\de;a_j,\infty;G_\pm)$$
    we obtain
    \begin{eqnarray}
    m(\Ga_4)&\le& m_x(\de,x-a_j;G_\pm)+C_{9}\nonumber\\
    \label{4.32}
    &\le&\frac{1}{2\pi}\log\frac{|J_j|}{\de}+C_{10}.
    \end{eqnarray}
    On the other hand, by virtue of (\ref{2.6}) and  (\ref{4.7})
    \beq\label{4.33}
    m(\phi_\pm(\Ga_4))\ge\frac{1}{\pi}\log\frac{t-a_j^*}{|t-\tilde{t}^\pm_\de|}\ge
    \frac{1}{\pi}\log\frac{|J_j|}{|t-\tilde{t}^\pm_\de|}-C_{11}.
    \eeq
    The inequalities (\ref{4.32}) and (\ref{4.33})
     imply (\ref{4.26}).

     \hfill$\Box$

  In the proof of Theorem \ref{th2} we need the following
  immediate consequence of Lemmas \ref{lem3.5}, \ref{lem3.6},
  \ref{lem4.7}, and  \ref{lem4.8}. Let $x\in E,0<\de\le 1,
  w=\phi(x), w_\de=\phi(x+i\de),
  t=\phi_\pm(x),\tilde{t}^\pm_\de=\phi_\pm(x\pm i\de)$. Then
  \beq\label{4.34}
  |w-w_\de|\ole|t-\tilde{t}^\pm_\de|.
  \eeq

 \absatz{The Extension Operator for an Antiderivative}

 Let $f\in BC^*_\om(E)$. We continuously extend $f$  to $\R$ such
 that on any component $J_j$ of $E^*$ it is a linear function.
 By the Lagrange formula for $x\in J_j$,
 $$|f(x)|=\left| f(b_j)\frac{x-a_j}{b_j-a_j}-
 f(a_j)\frac{x-b_j}{b_j-a_j}\right|\le||f||_{C(E)}.
 $$
 This clearly forces
 \beq\label{5.n1}
 ||f||_{C(\R)}=
 ||f||_{C(E)}.
 \eeq
 Furthermore, let $x_1,x_2\in E$ and $\xi_1,\xi_2\in\R$ satisfy
 $x_1\le\xi_1<\xi_2\le x_2$. For $k=1,2$ consider points $\nu_k\in
 E\cap [x_1,x_2]$ such that $|\xi_k-\nu_k|=d(\xi_k,E)$.
 Then
 \begin{eqnarray}
 |f(\xi_2)-f(\xi_1)|&\le& |f(\xi_2)- f(\nu_2)|+
 |f(\nu_2)-f(\nu_1)|+ |f(\nu_1)-f(\xi_1)|\nonumber\\
 \label{5.n2}
 &\le& 3\om(\mb{diam }\phi([x_1,x_2])).
 \end{eqnarray}
 Consider
 $$F(x):=\int_0^xf(t)dt\quad (x\in E).
 $$
 Our next objective is to continuously extend $F$  from $E$ to
 $\C$. The procedure described below is a modification of the
 corresponding constructions from \cite{ste} and \cite{dyn77}
 (see also  \cite[pp. 13-15]{andbeldzj}).
 \begin{lem}\label{lem5.1}
 (\cite[Chapter VI]{ste}). There exist a collection of closed
 squares $Q_k\subset\C\setminus E,k\in\N$ with sides
 parallel to the coordinate axes, and a set of infinitely
 differentiable (with respect to $x$ and $y$) functions
 $\mu_k(z)=\mu_k(x+iy)$ satisfying the following properties.

 (i) $\cup _kQ_k=\C\setminus E.$

 (ii) $\mb{\em diam }Q_k\le 2d(Q_k,E)\le 8 \mb{ \em diam }Q_k.$

 (iii) Each point $z\in\C\setminus E$ is contained in
  at most $144$ of the squares $ Q_k$.

 (iv) $\sum_k\mu_k(z)=1\quad (z\in\C\setminus E).$

 (v) $\mu_k(z)=0\quad (z\in\C\setminus Q_k).$

 (vi) For $r,l\in\N_0$ such that $r+l\le 1$ we have
 $$
 \left|\frac{\partial^{r+l}\mu_k(z)}{\partial x^r\partial
 y^l}\right|\ole (\mb{\em diam }Q_k)^{-r-l}.
 $$
 \end{lem}
 Let for $x\in E$ and $z\in\C$,
 $$P_x(z)=P_{x,f,E}(z):=F(x)+f(x)(z-x).
 $$
 Denote by $x_k\in E$ a point with the property
 $d(Q_k,E)=d(Q_k,x_k)$.
 The extension operator $\cal{E}$ is defined as follows:
 $$
 {\cal{E}}F(z):=
 \left\{\begin{array}{ll}
\sum_k'P_{x_k}(z)\mu_k(z) ,&\mb{ if }
z\in \C\setminus E,\\[2ex]
F(z),&\mb{ if }z\in E,
\end{array}\right.
 $$
where $\sum_k'$ means the sum taken  over only those squares $Q_k$
for which $d(Q_k,E)< 1$. Next, we clarify  the difference between
the sum $\Si_k'$ and the complete sum $\Si_k$. If $d(z,E)<1$, then
$\Si_k'=\Si_k$. If $d(z,E)>3$, then for every square $Q_k\ni z$,
by Lemma \ref{lem5.1}(ii),
$$3<d(z,E)\le\mb{ diam }Q_k+d(Q_k,E)\le 3d(Q_k,E),
$$
i.e., $\Si'_k$ does not have any terms at all.

For brevity, we also use the same notation $F$  for the extension
${\cal{E}}F$. The remark above yields
 \beq\label{5.2}
 F(z)=0\quad (z\in\C,\, d(z,E)\ge 3).
 \eeq
 \begin{lem}\label{lem5.2}
 The function $F$ possesses the following properties.

 (i) $F$ is continuous in $\C$.

 (ii) $F(z)=F(x+iy)$ is infinitely differentiable (with respect to
 $x$ and $y$) in $\C\setminus E$ and for $z\in\C\setminus E$ we have
 $$\left|\frac{\partial F(z)}{\partial\ov{z}}\right|\ole
  \left\{\begin{array}{ll}
\om(|\phi(z')-\phi(z'+id(z,E)|) ,&\mb{ if }
 d(z,E)<1,\\[2ex]
 |z|||f||_{C(E)},&\mb{ if }1\le d(z,E)\le 3,
\end{array}\right.
 $$
 where $z'\in E$ satisfies $d(z,E)=|z-z'|$.

 (iii) For $x\in E,z\in\C,|z-x|=\de,0<\de\le 1,$
 $$|F(z)-P_x(z)|\ole \om(|\phi(x)-\phi(x+i\de)|)\de.$$
 \end{lem}
 {\bf Proof.}
 Let $z\in Q_k$, by Lemma \ref{lem5.1}(iii) the number of such
 squares is at most $144$. According to Lemma \ref{lem5.1}(ii) for
 $\de_k:=$ diam $Q_k$ we obtain
 \beq\label{5.21}
 \frac{\de_k}{2}\le d(Q_k,E)\le d(z,E)\le \de_k+d(Q_k,E)\le
 5\de_k.
 \eeq
 Furthermore, if $z_k\in Q_k$ and $x_k\in E$ are such that
 $|x_k-z_k|=d(Q_k,E)$, then
 \beq\label{5.20}
 |z'-x_k|\le|z'-z|+|z-z_k|+|z_k-x_k|\le 11\de_k\le 22 d(z,E).
 \eeq

 (i) We only need to show that for any $x\in E$,
 \beq\label{5.4}
 \lim_{\C\setminus E\ni z\to x}F(z)=F(x).
 \eeq
 By virtue of Lemma \ref{lem5.1}(iv)
 for $z\in\C\setminus E$ with $|z-x|<1$ we have
 $$
 |F(z)-F(x)|=\sum_k(P_{x_k}(z)-F(x))\mu_k(z)
 $$
 which implies (\ref{5.4}).

 (ii) Consider the linear functions
 $$L_k(\z):=P_{x_k}(\z)-P_{z'}(\z)\quad (\z\in\C)
 $$
 and  points $u_k\in E$ such that
 $$|u_k-z'|\asymp|u_k-x_k|\asymp d(z,E).$$
 Since, by the Mean Value Theorem for $\xi_1,\xi_2\in E$,
 $$|P_{\xi_1}(\xi_2)-F(\xi_2)|= |f(\xi_1)-(\Re f(c_1)+i\Im
 f(c_2))||\xi_2-\xi_1|,
 $$
 for some $c_1$ and $c_2$ between $\xi_1$ and $\xi_2$,
 (\ref{1.5}), (\ref{3.6}), (\ref{5.n2}),(\ref{5.20}), and Lemma \ref{lem2.3}
  imply that
 \begin{eqnarray*}
 |L_k(x_k)|&=&|F(x_k)-P_{z'}(x_k)|\ole\om(\mb{diam }\phi([x_k,z']))|x_k-z'|\\
 &\ole& \om(|\phi(z')-\phi(z'+id(z,E)
 )|)d(z,E),
 \end{eqnarray*}
 \begin{eqnarray*}
 |L_k(u_k)|&\le&|P_{x_k}(u_k)-F(u_k|+|F(u_k)-P_{z'}(u_k)|\\
 &\ole& \om(\mb{diam }\phi([x_k,u_k]))|x_k-u_k|+
 \om(\mb{diam }\phi([u_k,z']))|u_k-z'|\\
 &\ole&\om(|\phi(z')-\phi(z'+id(z,E)
 )|)d(z,E),
 \end{eqnarray*}
 where $[x_k,z'],[x_k,u_k]$, and $[u_k,z']$ are the intervals of
 $\R$ between the appropriate points.

Moreover, by the Lagrange interpolation formula we obtain
 \begin{eqnarray}
 |L_k(z)|&\le& |L_k(u_k)|\left|\frac{z-x_k}{u_k-x_k}\right|+
 |L_k(x_k)|\left|\frac{z-u_k}{x_k-u_k}\right|\nonumber\\
 \label{5.3}
 &\ole& \om(|\phi(z')-\phi(z'+id(z,E)
 )|)d(z,E).
 \end{eqnarray}
 Therefore, if $d(z,E)<1$ then Lemma \ref{lem5.1}(vi),
 (\ref{5.21}), and (\ref{5.3}) yield that
 \begin{eqnarray}
 \left|\frac{\partial F(z)}{\partial\ov{z}}\right|&=&
 \left|\frac{\partial }{\partial\ov{z}}(F(z)-P_{z'}(z))\right|
 \nonumber\\
 \label{5.5}
 &=&\left|\sum_kL_k(z)\frac{\partial \mu_k(z)}{\partial\ov{z}}\right|
 \ole \om(|\phi(z')-\phi(z'+id(z,E)
 )|).
 \end{eqnarray}
 Let $1\le d(z,E)\le3$. Since according to (\ref{5.n1})
 $$|F(x)|\le |x|||f||_{C(E)}\quad (x\in E)
 $$
 and, therefore, for the squares $Q_k\ni z$,
 $$|P_{x_k}(z)|\ole|z|||f||_{C(E)},$$
 by virtue of Lemma \ref{lem5.1}(vi) we obtain
 $$\left|\frac{\partial F(z)}{\partial\ov{z}}\right|=
  \left|\sum'_kP_{x_k}(z) \frac{\partial \mu_k(z)}{\partial\ov{z}})\right|
  \ole
  |z|||f||_{C(E)},$$
  which, together with (\ref{5.5}), proves (ii).

  (iii) Since $F$ is continuous in $\C$ we can assume that $z\in C(x,\de)\setminus
  E$. For $Q_k\ni z$ consider the linear functions
  $$L_k^*(\z):=P_{x_k}(\z)-P_x(\z)\quad (\z\in\C)
  $$
  and  points $u_k^*\in E$ such that
  $$|u_k^*-x|\asymp|u_k^*-x_k|\asymp\de.
  $$
  Repeating the reasoning from the proof of the part (ii) we
  obtain
  $$|L_k^*(z)|\ole\om(|\phi(x)-\phi(x+i\de)|)\de.
  $$
  Therefore, according to  Lemma \ref{lem5.1}
  we have
  $$
  |F(z)-P_x(z)|=\left|\sum_kL_k^*(z)\mu_k(z)\right|\ole
  \om(|\phi(x)-\phi(x+i\de)|)\de,
  $$
  and (iii) is proved.

  \hfill$\Box$

\absatz{Auxiliary Entire Functions}

In this section we discuss the construction of certain entire
functions of exponential type. Our reasoning  is influenced by
\cite{shi03} and \cite[Chapter IX]{dzj}.

 We start with  well-known facts from the harmonic analysis
 (see, for details \cite{ach}). Let $g\in L(\R)$, i.e.,
 $$||g||_{L(\R)}:=\int_\R|g(x)|dx<\infty,$$
  and let
  $$\hat{g}(t):=\frac{1}{\sqrt{2\pi}}\int_\R g(x)e^{-itx}dx \quad
  (t\in\R)
  $$
  be the Fourier transform of $g$.

 For $\si>0$ and $s\in\N_0$ such that $0\le s<99$
 consider the function
 $$Q_{\si,s}(z):=\frac{\sin^{100}\si z}{z^{100-s}}\quad (z\in\C).$$
 Since $Q_{\si,s}\in E_{100\si}$ and
 $\int_\R|Q_{\si,s}(x)|^2dx<\infty$ we have
 $$\hat{Q}_{\si,s}(t)=0\quad (t\in\R,|t|\ge 100\si ).$$
 Moreover, for $t\in\R$,
  $$ g_{\si,s}(t):=\int_\R g(x)Q_{\si,s}(t-x)dx=
  \int_{-100\si}^{100\si}\hat{g}(x)\hat{Q}_{\si,s}(x)e^{itx}dx.
  $$
  Therefore, $g_{\si,s}$ can be extended  to the entire
  function (for which we use the same notation).
  Since
  $$||\hat{g}||_{C(\R)}\le\frac{1}{\sqrt{2\pi}}||g||_{L(\R)}$$
  and
  $$
 ||\hat{Q}_{\si,s}||_{C(\R)}\le\frac{1}{\sqrt{2\pi}}||Q_{\si,s}||_{L(\R)}\le
 C_1( \si),
 $$
 for $z\in\C$ we have
 \beq\label{6.1}
 |g_{\si,s}(z)|\le C_2(\si)||g||_{L(\R)}\exp(100 \si |\Im z|)
 \eeq
 (cf. \cite[p. 134]{ach}).

 Let $Q_{\si}:=Q_{\si,0}$. According to (\ref{6.1}) the function
 \begin{eqnarray*}
 g^*_{\si,s}(t)&:=&\int_\R g(x)x^sQ_{\si}(t-x)dx\\
 &=&\sum_{l=0}^s\frac{(-1)^ls!}{l!(s-l)!}t^{s-l}g_{\si,l}(t)
 \quad (t\in\R)
 \end{eqnarray*}
  belongs  to $E_{100\si
  }$. Moreover,  for $|z|\ge 1$,
  \beq\label{6.2}
 |g^*_{\si,s}(z)|\le C_3(\si)||g||_{L(\R)}|z|^s\exp(100 \si |\Im z|).
 \eeq

 Let $f\in BC^*_\om(E)$ and let $F$ be defined as in Section 5.
 For fixed $z_0:=4i$  let
 $$F_0(z):=\frac{F(z)}{(z-z_0)^3}\quad (z\in\C),$$
 $$\la_0(z):=\frac{1}{\pi}\frac{\partial
 F_0(z)}{\partial\ov{z}}=\frac{1}{\pi(z-z_0)^3} \frac{\partial
 F(z)}{\partial\ov{z}}\quad (z\in\C\setminus E),$$
 $$\la^\pm(w):=
 \left\{\begin{array}{ll}
 \la_0(\psi_\pm(w),&\mb{
if } w\in\He_\pm,
\psi_\pm(w)\in\He_\pm,\\[2ex]
 0,&\mb{ if } w\in\He_\pm,\psi_\pm(w)\in\He_\mp.
\end{array}\right.
$$
 Consider the kernel
 $$K(t):=C\left(\frac{\sin t}{t}\right)^{100}\quad (t\in\R),$$
 where
 $$C:=\left(\int_\R \left(\frac{\sin t}{t}\right)^{100}dt\right)^{-1}.$$
 For $\z\in G_{\pm},\si>0,$ and $t\in\R$ set
 $$\z^\pm_{\si,t}:=\psi_\pm\left(\phi_\pm(\z)\pm\frac{i}{\si}-t\right).$$
  Furthermore, for  $z\in \C\setminus G_\pm$ set
  $$
  e_\si^\pm(z):=
  \si\int_\R K(\si
  t)\int_{\He_\pm}\la_0(\z)\sum_{j=0}^{8}\frac{(\z_{\si,t}^\pm-\z)^{j}}
  {(\z_{\si,t}^\pm-z)^{j+1}}dm(\z) dt,
  $$
  where $dm(\z)$ means integration with respect to the
  two-dimensional Lebesgue measure (area).
   \begin{lem}\label{lem6.1}
 The function $e^\pm_{\si}$ can be extended  to the entire
 function belonging to $E_{C_4\si}$.
 \end{lem}
 {\bf Proof.} By virtue of (\ref{4.11}) and (\ref{5.2})
 for $z\in\C\setminus G_\pm$,
 \begin{eqnarray*}
 e_\si^\pm(z)&=& \si\int_\R K(\si t)\int_\R\int_0^{ v_0}\la^\pm(w)
 |\psi_\pm'(w)|^2\\
 &&\sum_{j=0}^{8}\frac{(\psi_\pm
 (w\pm\frac{i}{\si}-t)-\psi_\pm(w))^j}{ (\psi_\pm
 (w\pm\frac{i}{\si}-t)-z)^{j+1}}dvdudt,
 \end{eqnarray*}
 where $w=u\pm iv$ and $0<v_0\ole 1$.

 For $u\in\R,0<v<v_0$, and $0\le l\le j\le 8$ consider the functions
 $$\la_v^\pm(u):=\la^\pm(u\pm iv)|\psi_\pm'(u\pm iv)|^2,$$
 \begin{eqnarray*}
 &&h^\pm(z)=h^\pm_{\si,v,l,j}(z)\\
 &:=&
 \int_\R K(\si t)\int_\R\la^\pm_v(u)
 \frac{\psi_\pm
 (u\pm iv\pm\frac{i}{\si}-t)^l\psi_\pm(u\pm iv)^{j-l}}{ (\psi_\pm
 (u\pm iv\pm\frac{i}{\si}-t)-z)^{j+1}}dudt\\
 &=& \int_\R K(\si t)\int_\R\la^\pm_v(\xi+t)
 \frac{\psi_\pm
 (\xi\pm iv\pm\frac{i}{\si})^l\psi_\pm(t+\xi\pm iv)^{j-l}}{ (\psi_\pm
 (\xi\pm iv\pm\frac{i}{\si})-z)^{j+1}}d\xi dt\\
 &=& \int_\R\frac{\psi_\pm
 (\xi\pm iv\pm\frac{i}{\si})^l}{ (\psi_\pm
 (\xi\pm iv\pm\frac{i}{\si})-z)^{j+1}}q^\pm(\xi) d\xi,
 \end{eqnarray*}
 where
 \begin{eqnarray*}
 q^\pm(\xi)&:=&\int_\R K(\si t)p^\pm(\xi+t)dt\\
 &=& \int_\R K(\si (\eta-\xi))r^\pm(\eta)(\eta\pm i)^8
 d\eta,
 \end{eqnarray*}
 \begin{eqnarray*}
 p^\pm(\xi)&:=&\la^\pm_v(\xi)\psi_\pm(\xi\pm iv)^{j-l}\\
 &=&\la^\pm(\xi\pm iv)|\psi_\pm'(\xi\pm iv)|^2 \psi_\pm(\xi\pm
 iv)^{j-l},
 \end{eqnarray*}
 $$
 r^\pm(\eta):=\frac{p^\pm(\eta)}{(\eta\pm i)^8}.
 $$
 By (\ref{4.1}) and Lemma \ref{lem5.2}(ii)
 $$|\la^\pm(\xi\pm iv)|\ole ||f||_{C(E)}|\xi-z_0|^{-2},$$
 and according to (\ref{4.12}) for $\xi\pm iv$ such that $\la^\pm(
 \xi\pm iv)\neq 0$ we obtain
 $$|\psi_\pm'(\xi\pm iv)|\ole 1.
 $$
 Therefore, by virtue of (\ref{4.1}) and (\ref{6.2}) for $|\z|\ge1$
 we have
 \beq\label{6.3}
 |q^\pm(\z)|\le
 C_5(E,\si)||f||_{C(E)}|\z|^{8}
 \exp(100\si |\Im\z|).
 \eeq
 Notice that $h^\pm$ can be extended analytically from
 $\C\setminus G_\pm$ to $\C$ as follows. By (\ref{2.14}) for $z\in \C$ with
 $|z|\ge 1$ we can find $r$ such that
 $$2|z|\le|\psi_\pm(w)|\ole|z|\quad (w\in \ov{\He_\pm}\cap C(r)).$$
 According to (\ref{4.1}) $r\le C_6|z|$.
 Then the above mentioned extension is defined by the formula
 $$h^\pm(z)=\int_{J_r}\frac{\psi_\pm(\xi\pm iv
 \pm\frac{i}{\si})^l}{(\psi_\pm(\xi\pm iv
 \pm\frac{i}{\si})-z)^{j+1}}q^\pm(\xi)d\xi,
 $$
 where $J_r:=(C(r)\cap\He_\pm)\cup(\R\setminus D(r))$.
 We may now differentiate $9$ times to conclude that
 according  to
 (\ref{4.1}) and (\ref{6.3}) for $0<v<v_0$,
 $$\left|\frac{d^9}{dz^9}h^\pm(z)\right|\ole\exp(100C_6\si|z|).$$
 Therefore,
 $$
 |h^\pm(z)|\ole\exp(C_7\si|z|).
 $$
 Comparing the definitions of $e_\si^\pm$ and
 $h^\pm$ with the above inequality we have
 $e^\pm_\si\in E_{C_7\si}.$

 \hfill$\Box$

\absatz{Proof of Theorem \ref{th2}}

Let $z_0=4i$ and let functions
$F,P_x,\la_0,e_\si^\pm,e_\si:=e_\si^++e_\si^-$ be defined as in
Sections 5 and 6. By (\ref{1.5}) and Lemma \ref{lem6.1} in order
to prove (\ref{1.6}) it is sufficient to show that  for $\si>2$
and $x\in E$,
 \beq\label{7.1}
 \left|\frac{f(x)}{(x-z_0)^3}-e_\si'(x)-\frac{3e_\si(x)}{x-z_0}\right|\ole
 |x-z_0|^{-3}\left(\frac{||f||_{C(E)}}{\si}+\om\left(\frac{1}{\si}\right)\right).
 \eeq
 Let for $\z\in G_\pm $ and $x\in E$,
 $$R^\pm(\z,x)=R^\pm_{\si,t}(\z,x):=\sum_{j=0}^{8}
  \frac{(\z^\pm_{\si,t}-\z)^j}{(\z^\pm_{\si,t}-x)^{j+1}},
 $$
 $$
 R(\z,x)=R_{\si,t}(\z,x):=\sum_\pm R^\pm(\z,x),
 $$
 $$
 d:=d_\si(x),\quad D:=D(x,d),\quad J:=C(x,d),\quad
 \Si:=\{z:\,d(z,E)\le 3\}.
 $$
 Since for $x\in E$
 \begin{eqnarray*}
 \frac{f(x)}{(x-z_0)^3}&=&\frac{1}{2\pi i}\int_J
 \frac{F(\z)}{(\z-z_0)^3}\left(\frac{1}{(\z-x)^2}+\frac{3}{(\z-x)(x-z_0)}\right)
 d\z\\
 &+& \frac{1}{2\pi i}\int_J
 \frac{P_x(\z)-F(\z)}{(\z-z_0)^3}\left(\frac{1}{(\z-x)^2}+\frac{3}{(\z-x)
 (x-z_0)}\right) d\z
 \end{eqnarray*}
 and according to the Green formula (cf. \cite[p. 22]{andbeldzj})
 \begin{eqnarray*}
 &&\frac{1}{2\pi i}\int_J
 \frac{F(\z)}{(\z-z_0)^3}\left(\frac{1}{(\z-x)^2}+\frac{3}{(\z-x)(x-z_0)}\right)
 d\z\\
 &=&
 \int_{\Si\setminus D}\la_0(\z)\left(\frac{\partial}{\partial
 x}\frac{1}{\z-x}+\frac{3}{x-z_0}\frac{1}{\z-x}\right) dm(\z),
 \end{eqnarray*}
 we have
 \begin{eqnarray}
 &&\frac{f(x)}{(x-z_0)^3}-e_\si'(x)-\frac{3e_\si(x)}{x-z_0}
 \nonumber\\
 \label{7.2}
 &=& \si\int_\R K(\si t)\left(\sum_\pm\sum_{l=0}^4I^\pm_l\right)
 dt+I_5,
 \end{eqnarray}
 where
 $$\Si_0:=\{z\in\C\setminus E:\, d(z,E)<1\},$$
 \begin{eqnarray*}
 I^\pm_0&:=&\int_{\He_\pm\cap\Si\setminus\Si_0}\la_0(\z)
 \left[\frac{\partial}{\partial
 x}\left(\frac{1}{\z-x}-R^\pm(\z,x)\right)\right.\\
 &&+\frac{3}{x-z_0}\left. \left(
 \frac{1}{\z-x}-R^\pm(\z,x)\right)\right]dm(\z),\\
 I^\pm_1&:=&\int_{\He_\pm\cap\Si_0\setminus D}\la_0(\z)
 \frac{\partial}{\partial
 x}\left(\frac{1}{\z-x}-R^\pm(\z,x)\right)dm(\z),\\
 I^\pm_2&:=&-\int_{\He_\pm\cap D}\la_0(\z)
 \frac{\partial}{\partial
 x} R^\pm(\z,x) dm(\z),\\
 I^\pm_3&:=&\frac{3}{x-z_0}\int_{\He_\pm\cap D}\la_0(\z)
 \left(\frac{1}{\z-x}-R^\pm(\z,x)\right)dm(\z),\\
 I^\pm_4&:=&-\frac{3}{x-z_0}\int_{\He_\pm D}\la_0(\z)
 R^\pm(\z,x)dm(\z),
 \end{eqnarray*}
 $$I_5:=\frac{1}{2\pi
 i}\int_J\frac{P_x(\z)-F(\z)}{(\z-z_0)^3}\left(\frac{1}{(\z-x)^2}
 +\frac{3}{(\z-x)(x-z_0)}\right) d\z.
 $$
  According to Lemma \ref{lembel} for $\z\in G_\pm$ and $x\in E$,
  \begin{eqnarray}
  \left|\frac{\z^\pm_{\si,t}-\z}{\z^\pm_{\si,t}
  -x}\right|&=&
   \left|\frac{\z^\pm_{\si}-\z}{\z^\pm_{\si}
  -x}\right|\left|\frac{\z^\pm_{\si,t}-\z}{\z^\pm_{\si}
  -\z}\right|\left|\frac{\z^\pm_{\si}-x}{\z^\pm_{\si,t}
  -x}\right|\nonumber\\
  \label{7.3}
  &\ole& (1+\si |t|)^8\left|\frac{\z^\pm_{\si}-\z}{\z^\pm_{\si}
  -x}\right|,
  \end{eqnarray}
  \beq\label{7.3dob}
  \left|\frac{\z^\pm_{\si}-x}{\z^\pm_{\si,t}
  -x}\right|\ole (1+\si |t|)^4.
  \eeq
 Since  by Lemma \ref{lem2.3} for $\z\in\He_\pm$ and $x\in E$,
 $$|\z-x|\ole |\z^\pm_{\si}-x|$$
 and
 $$
 \frac{1}{\z-x}-R^\pm(\z,x)=\frac{1}{\z-x}
 \left(\frac{\z^\pm_{\si,t}-\z}{\z^\pm_{\si,t}-x}\right)^9,
 $$
 \begin{eqnarray*}
 &&\frac{\partial}{\partial x}\left[
 \frac{1}{\z-x}-R^\pm(\z,x)\right]\\
 &=&\frac{1}{(\z-x)^2}
 \left(\frac{\z^\pm_{\si,t}-\z}{\z^\pm_{\si,t}-x}\right)^9+
 \frac{9}{\z-x}
 \frac{(\z^\pm_{\si,t}-\z)^9}{(\z^\pm_{\si,t}-x)^{10}},
 \end{eqnarray*}
  according to  (\ref{7.3}), (\ref{7.3dob}), and Lemma \ref{lem4.5} for
   $\z\in\He_\pm\cap\Si\setminus D$ we have
  \beq\label{7.4}
  \left|
 \frac{1}{\z-x}-R^\pm(\z,x)\right|\ole\frac{1}{|\z-x|}
 \left(\frac{d}{|\z -x|}\right)^{9/4}(1+\si |t|)^{72},
 \eeq
  \beq\label{7.5}
  \left| \frac{\partial}{\partial x}\left(
 \frac{1}{\z-x}-R^\pm(\z,x)\right)\right|\ole\frac{1}{|\z-x|^2}
 \left(\frac{d}{|\z -x|}\right)^{9/4}(1+\si |t|)^{76}.
 \eeq
 Moreover, by virtue of Lemmas \ref{lembel} and \ref{lem4.4} for
 $\z\in\He_\pm\cap D$,
  \beq\label{7.6}
  \left|
 R^\pm(\z,x)\right|\ole (1+\si |t|)^{68}
 \sum_{j=1}^{8}
 \frac{|\z_\si^\pm-\z|^j}{|\z_\si^\pm-x|^{j+1}}
 \ole \frac{(1+\si|t|)^{68}}{d},
 \eeq
  \beq\label{7.7}
  \left|\frac{\partial}{\partial x}
 R^\pm(\z,x)\right|\ole (1+\si |t|)^{72}
 \sum_{j=1}^{8}
 \frac{|\z_\si^\pm-\z|^j}{|\z_\si^\pm-x|^{j+2}}
 \ole \frac{(1+\si|t|)^{72}}{d^2}.
 \eeq
 According to (\ref{1.5}), (\ref{4.34}), and Lemma \ref{lem5.2}(ii) we obtain
 $$
 |\la_0(\z)||\z-z_0|^3\ole
 $$
 \beq\label{7.8}
 \left\{\begin{array}{ll}
\om(|\phi_\pm(x)-\phi_\pm(x\pm i|x-\z|)|),&\mb{ if }
 \z\in \He_\pm\cap\Si_0,\\[2ex]
 |\z|||f||_{C(E)},&\mb{ if }\z\in \He_\pm\cap\Si\setminus\Si_0.
\end{array}\right.
\eeq
 Next, applying (\ref{7.4})-(\ref{7.8}) we estimate each of
 the integrals
 $$
 \tilde{I}_l^\pm:=\left|\sigma\int_\R K(\si t)I_l^\pm
 dt\right|\quad
 (l=0,\ldots,4)
 $$
  from above.

 Since by (\ref{4.231}) and Lemma \ref{lem4.7}
 \beq\label{7.88}
 d\ole \si^{-1}
 \eeq
 and
 $$\tilde{I}_0^\pm\ole\frac{||f||_{C(E)}}{\si^{9/4}}\int_{
 \He_\pm\cap\Si\setminus\Si_0}\left(\frac{1}{|\z-x|}+ \frac{1}{|x-z_0|}
 \right) \frac{dm(\z)}{|\z-z_0|^2|\z-x|^{13/4}},
 $$
 dividing the set of integration into subsets
 \begin{eqnarray*}
 U_1^\pm&:=& \He_\pm\cap(\Si\setminus\Si_0)\cap
 D\left(x,\frac{1}{2}|x-z_0|\right),\\
 U_2^\pm&:=& \He_\pm\cap(\Si\setminus\Si_0)\cap
 D\left(x,2|x-z_0|\right)\setminus U_1^\pm,\\
 U_3^\pm&:=& \He_\pm\cap(\Si\setminus\Si_0) \setminus(U_1^\pm\cup
 U_2^\pm),
 \end{eqnarray*}
 and passing to the polar coordinates with the center either in $x$
 or in $z_0$ we obtain
 \beq\label{7.9}
 \tilde{I}^\pm_0\ole
 \frac{||f||_{C(E)}}{\si |x-z_0|^3}.
 \eeq
 To deal with $\tilde{I}^\pm_1$, we note that by virtue of (\ref{1.5}),
 (\ref{4.23}), and
 (\ref{4.201})
 $$
 \frac{\om(|\phi_\pm(x)-\phi_\pm(x\pm i|x-\z|)|)}{\om
 (|\phi_\pm(x)-\phi_\pm(x\pm id)|)}\ole\left|
 \frac{\phi_\pm(x)-\phi_\pm(x\pm i|x-\z|)}{
 \phi_\pm(x)-\phi_\pm(x\pm id)}\right|\ole\frac{|x-\z|}{d}.
 $$
 Therefore, according to (\ref{1.5}) and (\ref{2.14})
 $$
 \tilde{I}^\pm_1\ole\om\left(\frac{1}{\si}\right)
 d^{5/4}\int_{\He_\pm\cap\Si_0\setminus D}\frac{dm(\z)}{|\z-z_0|^3
 |\z-x|^{13/4}}.
 $$
 Next, dividing the set of integration $\He_\pm\cap\Si_0\setminus
 D$ into subsets
  \begin{eqnarray*}
 V_1^\pm&:=& \He_\pm\cap(\Si_0\setminus D)\cap
 D\left(x,\frac{1}{2}|x-z_0|\right),\\
 V_2^\pm&:=& \He_\pm\cap(\Si_0\setminus D)\cap
 D\left(x,2|x-z_0|\right)\setminus V_1^\pm,\\
 V_3^\pm&:=& \He_\pm\cap(\Si_0\setminus D) \setminus(V_1^\pm\cup
 V_2^\pm)
 \end{eqnarray*}
 and passing to the polar coordinates with the center either in $x$
 or in $z_0$ as well as applying (\ref{7.88}) we obtain
 \beq\label{7.10}
 \tilde{I}^\pm_1\ole \om\left(\frac{1}{\si}\right)|x-z_0|^{-3}.
 \eeq
 In the same manner we can prove that
 \beq\label{7.11}
 \tilde{I}^\pm_3\ole \om\left(\frac{1}{\si}\right)|x-z_0|^{-3}.
 \eeq
 In order to estimate $\tilde{I}^\pm_2$ and $\tilde{I}^\pm_4$ we
 note that by (\ref{2.14}) and (\ref{4.16}) for $\z\in\He_\pm\cap D$,
 $$
 |\phi_\pm(x)-\phi_\pm(x\pm i|x-\z|)|\ole\frac{1}{\si},
 $$
 $$
 |\z^\pm_\si-\z|\oge d.
 $$
 Therefore,
 \beq\label{7.12}
 \tilde{I}^\pm_2\ole \om\left(\frac{1}{\si}\right)|x-z_0|^{-3}
 d^{-2}\int_{\He_\pm\cap D}dm(\z)\ole
 \om\left(\frac{1}{\si}\right)|x-z_0|^{-3},
 \eeq
 \beq\label{7.13}
 \tilde{I}^\pm_4\ole \om\left(\frac{1}{\si}\right)|x-z_0|^{-4}
 d^{-1}\int_{\He_\pm\cap D}dm(\z)\ole
 \om\left(\frac{1}{\si}\right)|x-z_0|^{-3}.
 \eeq
 Furthermore,  (\ref{4.34}) and Lemma \ref{lem5.2}(iii) imply that
 \begin{eqnarray}
 |I_5|&\ole&\om\left(\frac{1}{\si}\right)|x-z_0|^{-3} d
 \int_J\left(\frac{1}{|\z-x|^2}+\frac{1}{|\z-x||x-z_0|}\right)
 |d\z|\nonumber\\
 \label{7.14}
 &\ole&\om\left(\frac{1}{\si}\right)|x-z_0|^{-3}.
 \end{eqnarray}
 Comparing (\ref{7.2}) and (\ref{7.9})-(\ref{7.14}) we have
 (\ref{7.1}).

 \hfill$\Box$

\absatz{Proof of Theorem \ref{th3}}

 We adapt to our case the standard procedure of proving  inverse theorems.
  Let $x_1,x_2\in E$ be such that $x_1<x_2$ and let
 $\de:=\rho_E(x_1,x_2)$.
  Since
 $$
 |f(x_2)-f(x_1)|\le 2
 ||f||_{C(E)},
 $$
 we can assume that $0<\de<1/2$.

 Let $e_k\in E_{2^k},k\in\N_0$ satisfy
 $$||f-e_k||_{C(E)}\le 2\om(2^{-k})
 $$
  and let
 $$
 g_k(z):=e_{k+1}(z)-e_{k}(z).
 $$
 Chose $m\in\N$ such that
 $$2^{-m-1}<\de\le 2^{-m}$$
 Since
 \begin{eqnarray*}
 f(x_2)-f(x_1)&=&[f(x_2)-e_m(x_2)]-[f(x_1)-e_m(x_1)]\\
 &+& e_0(x_2)-e_0(x_1)+\sum_{k=0}^{m-1}[g_k(x_2)
 -g_k(x_1)],
 \end{eqnarray*}
 we have
 \beq\label{8.2}
 |f(x_2)-f(x_1)|\le 4\om(2\de)+\int_{x_1}^{x_2}|e_0'(x)|dx
 +\sum_{k=0}^{m-1} \int_{x_1}^{x_2}|g_k'(x)|dx.
 \eeq
 Furthermore, since for $x\in E$ and $k=0,\ldots,m-1$,
 $$
 |e_0(x)|\le |f(x)|+ |f(x)-e_0(x)|\le ||f||_{C(E)}+2\om(1),
 $$
 $$
 |g_k(x)|\le|e_{k+1}(x)-f(x)|+|f(x)-e_{k}(x)|\le 4\om(2^{-k}),
 $$
 by the Phragmen-Lindel\"of theorem (see \cite[VIIIA]{koo})
 for $z\in\C$,
 \beq\label{8.3}
 |e_0(z)|\le (||f||_{C(E)}+2\om(1)) \exp(C_1|\Im\phi(z)|),
 \eeq
  \beq\label{8.4}
  |g_k(z)|\le 4\om(2^{-k})
   \exp(C_12^{k+1}|\Im\phi(z)|),
 \eeq
 where $\phi$ is the Levin conformal mapping extended to $\He_-$
 by the formula
 $$
 \phi(z):=\ov{\phi(\ov{z})}\quad (z\in\He_-).
 $$
 Let for $x\in E$,
 $$E^\pm_{\de}:=\{z\in\ov{\He}_\pm:\, \pm \Im\phi(z)=\de),\quad
 \rho_{\de}(x)= d(x,E^\pm_{\de}).
 $$
 According to  (\ref{3.20}) for $x\in E$ and
 $0<t<T$,
 \beq\label{8.5}
 \left(\frac{t}{T}\right)^{C_2}\ole\frac{\rho_{t}(x)}{\rho_{T}(x)}
 \ole \frac{t}{T}.
 \eeq
 Since by (\ref{1.31})
 $$
 x_2-x_1<C_3\de<\frac{C_3}{2}
 $$
 and according to (\ref{3.19}) and (\ref{8.5})
 $$2(x_2-x_1)\le\rho_{C\de}(x_1),
 $$
 by virtue of (\ref{1.5}), (\ref{3.2}), and (\ref{8.3})-(\ref{8.5})
  for $x_1\le x\le x_2$ we have
 $$
 |e_0'(x)|\le\frac{1}{2\pi}\int
 _{|z-x_1|=C_3}\frac{|e_0(z)|}{|z-x|^2}|dz| \ole
 ||f||_{C(E)}+\frac{\om(\de)}{\de},
 $$
 $$
 |g_k'(x)|\le\frac{1}{2\pi}\int
 _{|z-x_1|=\rho_{C2^{-k}}(x_1)}\frac{|g_k(z)|}{|z-x|^2}|dz| \ole
  \frac{\om(2^{-k})}{\rho_{2^{-k}}(x_1)}.
 $$
 Therefore,
 \beq\label{8.8}
 \int_{x_1}^{x_2}|e_0'(x)|dx\ole \de ||f||_{C(E)}+\om(\de),
 \eeq
 \begin{eqnarray}
 \sum_{k=0}^{m-1}\int_{x_1}^{x_2}|g_k'(x)|dx&\ole&
  \sum_{k=0}^{m-1}\frac{\om(2^{-k})\rho_{2^{-m}}(x_1)}{
  \rho_{2^{-k}}(x_1)}\nonumber\\
  \label{8.9}
 &\ole& \sum_{k=0}^{m-1}2^{k-m}\om(2^{-k})\ole\de
 \int_\de^1\frac{\om(t)}{t^2}dt.
 \end{eqnarray}
 Comparing (\ref{8.2}), (\ref{8.8}), and (\ref{8.9}) we obtain
 (\ref{1.8}).

 \hfill$\Box$

 \absatz{Acknowledgements}

 This research was conducted  while the author was visiting the Katholische
Universit\"at Eichstaett and Universit\"at Trier. The author
wishes to thank the members of these universities for the pleasant
mathematical atmosphere they offered him.  The author would also
like to warmly thank  Misha Nesterenko for  many useful remarks.

V. V. Andrievskii

 Department of Mathematical Sciences

 Kent State University

 Kent, OH 44242

 USA

e-mail: andriyev@math.kent.edu

\end{document}